\documentclass[leqno]{amsart}
\usepackage{amsmath,amssymb,xspace,upref,url}

\usepackage[all,cmtip,ps]{xy}
\newdir^{ (}{{}*!/-5pt/\dir^{(}}
\newdir_{ (}{{}*!/-5pt/\dir_{(}}
\let\hb z
\let\hbm\zeta
\let\leq\leqslant
\let\geq\geqslant
\let\wt\widetilde
\let\ov\overline

\let\moins\smallsetminus
\let\epsilon\varepsilon
\def\eg{e.g.,\ }
\def\ie{i.e.,\ }
\def\cf{cf.\kern.3em}
\def\resp{resp.\kern.3em}
\def\loccit{loc.\kern2pt cit.\xspace}
\def\defin{:=}
\def\onabla{{}^\omega\!\nabla}

\let\oldboxtimes\boxtimes
\def\sboxtimes{{\scriptstyle\oldboxtimes}}
\def\ssboxtimes{{\scriptscriptstyle\oldboxtimes}}
\def\boxtimes{\mathchoice{\mathrel{\sboxtimes}}{\mathrel{\sboxtimes}}{\ssboxtimes}{\ssboxtimes}}

\def\oR{{}^{\otimes r\!}R}
\def\Ur{U^r\!}
\def\gE{\mathfrak{E}}
\def\gS{\mathfrak{S}}
\def\gm{\mathfrak{m}}

\def\bun{\mathbf{1}}
\def\AA{\mathbb{A}}
\def\CC{\mathbb{C}}
\def\NN{\mathbb{N}}
\def\PP{\mathbb{P}}
\def\QQ{\mathbb{Q}}
\def\ZZ{\mathbb{Z}}
\def\cD{\mathcal{D}}
\def\cF{\mathcal{F}}
\def\cG{\mathcal{G}}
\def\cH{\mathcal{H}}
\def\cL{\mathcal{L}}
\def\cM{\mathcal{M}}

\def\cO{\mathcal{O}}
\def\cU{\mathcal{U}}

\def\tot{\mathrm{tot}}
\def\rw{\mathrm{w}}
\def\rW{\mathrm{W}}

\def\nablab{\boldsymbol{\nabla}}

\def\bme{\boldsymbol{e}}
\def\bE{\boldsymbol{E}}
\def\bF{\boldsymbol{F}}
\def\bR{\boldsymbol{R}}
\newcommand{\cbbullet}{{\raisebox{1pt}{$\scriptscriptstyle\bullet$}}}
\let\oldvee\vee
\def\vee{{\scriptscriptstyle\oldvee}}
\let\oldforall\forall
\def\forall{\oldforall\,}
\let\oldexists\exists
\def\exists{\oldexists\,}
\def\isom{\stackrel{\sim}{\longrightarrow}}

\def\Afu{\AA^{\!1}}
\def\pQQ{{{}^p\!\QQ}}
\def\pCC{{{}^p\!\CC}}
\def\pcH{{{}^p\!\cH}}
\def\phip{{{}^p\!\phi}}

\DeclareMathOperator{\ad}{ad}
\DeclareMathOperator{\coker}{coker}
\DeclareMathOperator{\diag}{diag}

\DeclareMathOperator{\Hom}{Hom}

\DeclareMathOperator{\id}{Id}

\DeclareMathOperator{\im}{Im}
\DeclareMathOperator{\Perv}{Perv}
\DeclareMathOperator{\rk}{rk}
\DeclareMathOperator{\sgn}{sgn}

\DeclareMathOperator{\Sp}{Sp}

\DeclareMathOperator{\vol}{vol}
\newcommand{\ant}{\mathrm{ant}}
\newcommand{\Elem}{\mathrm{Elem}}
\newcommand{\Sym}{\mathrm{Sym}}

\DeclareMathAlphabet{\othercal}{U}{eus}{m}{n}

\def\to{\mathchoice{\longrightarrow}{\rightarrow}{\rightarrow}{\rightarrow}}
\def\mto{\mathchoice{\longmapsto}{\mapsto}{\mapsto}{\mapsto}}
\def\hto{\mathrel{\lhook\joinrel\to}}
\def\To#1{\mathchoice{\xrightarrow{\textstyle\kern4pt#1\kern3pt}}{\stackrel{#1}{\longrightarrow}}{}{}}
\def\Hto#1{\mathrel{\lhook\joinrel\To{#1}}}

\numberwithin{equation}{section}
\theoremstyle{plain}
\newtheorem{lemma}[equation]{Lemma}
\newtheorem{proposition}[equation]{Proposition}
\newtheorem{corollary}[equation]{Corollary}
\newtheorem{theorem}[equation]{Theorem}
\newtheorem{assumption}[equation]{Assumption}

\theoremstyle{definition}
\newtheorem{definition}[equation]{Definition}

\theoremstyle{remark}
\newtheorem{remark}[equation]{Remark}
\newtheorem{remarks}[equation]{Remarks}
\newtheorem{example}[equation]{Example}

\newenvironment{enumeratei}
{\bgroup\begin{enumerate}}
{\end{enumerate}\egroup}
\newenvironment{enumeratea}
{\bgroup\begin{enumerate}}
{\end{enumerate}\egroup}

\begin{document}
\title[The Grassmannian and alternate Thom-Sebastiani]{Quantum cohomology of the Grassmannian and alternate~Thom-Sebastiani}

\author{Bumsig Kim}
\email{bumsig@kias.re.kr}
\address{School of Mathematics\\
Korea Institute for Advanced Study\\
207-43 Cheongnyangni 2-dong, Dongdaemun-gu\\
Seoul, 130-722\\
Korea}

\author{Claude Sabbah}
\email{sabbah@math.polytechnique.fr}
\address{UMR 7640 du C.N.R.S.\\
Centre de math\'ematiques Laurent Schwartz,
\'Ecole polytechnique\\
F-91128 Palaiseau cedex\\
France}



\begin{abstract}
We introduce the notion of alternate product of Frobenius manifolds and we give, after \cite{CF-K-S06}, an interpretation of the Frobenius manifold structure canonically attached to the quantum cohomology of $G(r,n+1)$ in terms of alternate products. We also investigate the relationship with the alternate Thom-Sebastiani product of Laurent polynomials.
\end{abstract}

\thanks{The first author thanks the staff at \'Ecole polytechnique for the warm hospitality during his visit. His work is supported by KOSEF grant R01-2004-000-10870-0.  The second author thanks KIAS for providing him with excellent working conditions during his visit.}

\maketitle

\section*{Introduction}
It is known that the Frobenius manifold structure attached canonically to the quantum cohomology of the complex projective space $\PP^n$ can also be obtained, in a canonical way, by considering the Laurent polynomial $f(u_1,\dots,u_n)=u_1+\cdots+u_n+1/u_1\cdots u_n$ on the torus $U=(\CC^*)^n$ and its associated Gauss-Manin system (\cf \cite{Barannikov00}).

The main result of \cite{CF-K-S06} applied to the case of the complex Grassmann variety $G(r,n+1)$ of $r$-planes in $\CC^{n+1}$ explains how to compute the Frobenius manifold structure canonically attached to the quantum cohomology of $G(r,n+1)$ in terms of that of $\PP^n$.

In this article, we introduce the notion of alternate product of Frobenius manifolds and we give an interpretation of the previous result in terms of alternate products.

On the ``mirror side'', let us consider the following data:
\begin{itemize}
\item
the affine variety $U^{(r)}$ obtained as the quotient of the $r$-fold product $U^r$ by the symmetric group $\gS_r$,
\item
the function $f^{(\oplus r)}$ on $U^{(r)}$ induced by the $r$-fold Thom-Sebastiani sum $f^{\oplus r}:U^r\to\CC$,
\item
the rank-one local system $\cL$ on the complement of the discriminant (image of the diagonals) in $U^{(r)}$, corresponding to the signature $\sgn:\gS_r\to\{\pm1\}$.
\end{itemize}
We show that the Gauss-Manin system attached to these data is the $r$-fold alternate product of that of $f$, making these data a candidate for being ``mirror of the Grassmannian''.

The contents of the article is as follows: Section \ref{sec:1} recalls the correspondence between Frobenius and Saito structures on a manifold. The point of view of Saito structures (primitive forms) enables us to use the results of Hertling and Manin \cite{H-M04} to generate Frobenius manifold structures.

This construction is applied to tensor and alternate products in Section \ref{sec:tensor}. We express the quantum cohomology of the Grassmannian, as a Frobenius manifold, in terms of the alternate product of that of the projective space in Theorem \ref{th:cfks}, which is mainly a reformulation of \cite[Th.~4.1.1(a)]{CF-K-S06} in this context.

In Section \ref{sec:TS}, we show that the Gauss-Manin system, with coefficients in the local system $\cL$, of the function $f^{(r)}$ considered above, can be obtained as the $r$-fold alternate product of the Gauss-Manin system of $f$.

In Section \ref{sec:geom}, we recall the notion of canonical Frobenius manifold attached to a Laurent polynomial satisfying generic assumptions, and we conclude that the Gauss-Manin system of the pair $(f^{(\oplus r)},\cL)$ on $U^{(r)}$ is also obtained from the quantum cohomology of the Grassmannian.

\subsubsection*{Acknowledgements}
The authors thank the referee for his careful reading of the manuscript.

\section{Saito and Frobenius manifold structures}\label{sec:1}
In this section, we will work in the category of punctual germs of complex analytic manifolds, although most of the results can be extended to simply connected complex analytic manifolds. We denote by $\cO_M$ the local algebra of $M$, by $\gm$ its maximal ideal and by $\Theta_M$ the tangent bundle of $M$.

\subsection{Pre-Saito structures}\label{subsec:prefro}
We refer to \cite[\S VI.2.c]{Bibi00} for more details on what follows. By a pre-Saito structure (without metric) on $M$ we mean a t-uple $(M,E,\nabla,R_\infty,\Phi,R_0)$ where
\begin{itemize}
\item
$E$ is a vector bundle on $M$,
\item
$\nabla$ is a connection on $E$,
\item
$R_0,R_\infty$ are $\cO_M$-linear endomorphisms of $E$,
\item
$\Phi:\Theta_M\otimes_{\cO_M}E\to E$ is a $\cO_M$-linear morphism,
\end{itemize}
which satisfy the following relations:
\begin{gather*}
\nabla^2=0,\quad
\nabla(R_\infty)=0,\quad
\Phi\wedge\Phi=0,\quad
{[R_0,\Phi]}=0\\
\nabla(\Phi)=0,\quad
\nabla(R_0)+\Phi=[\Phi,R_\infty].
\end{gather*}
In particular, $\nabla$ is flat and $\Phi$ is a Higgs field. These conditions are better understood by working on the manifold $M\times\Afu$, where $\Afu$ is the affine line with coordinate $\hb$. Let $\pi:M\times\Afu\to M$ denote the projection. Then, on $\bE\defin\pi^*E$, the connection $\nablab$ defined by
\begin{equation}\label{eq:nablab}
\nablab=\pi^*\nabla+\hb\Phi+(R_\infty-\hb R_0)\frac{d\hb}{\hb}
\end{equation}
is flat if and only if the previous relations are satisfied. We will also denote a pre-Saito structure by $(\bE,\nablab)$.

Let us fix local coordinates $x_1,\dots,x_m$ on $M$ and let $\bme$ be a $\nabla$-horizontal basis of $E$. We also set\footnote{The use of $-B_\infty$ instead of $B_\infty$ is done to keep a perfect correspondence with \cite[Chap.~VI]{Bibi00}.}
\[
R_\infty(\bme)=\bme\cdot(-B_\infty),\quad\Phi_{\partial_{x_i}}(\bme)=\bme\cdot C^{(i)}(x),\quad R_0(\bme) =\bme\cdot B_0(x).
\]
Then the previous relations reduce to the constancy of $B_\infty$ and to
\begin{equation}\label{eq:relprefrobhoriz}
\begin{aligned}
\frac{\partial C^{(i)}}{\partial x_j}&=\frac{\partial C^{(j)}}{\partial x_i},\\
[C^{(i)},C^{(j)}]&=0,\\
[B_0,C^{(i)}]&=0\\
C^{(i)}+\frac{\partial B_0}{\partial x_i}&=[B_\infty,C^{(i)}].
\end{aligned}
\end{equation}

\subsection{Universal deformation}\label{subsec:univdef}
Let $f:N\to M$ be a holomorphic map and let $Tf:\Theta_N\to f^*\Theta_M$ be its tangent map. Then the pull-back of a t-uple $(M,E,\nabla,R_\infty,\Phi,R_0)$ is defined by
\begin{itemize}
\item
$f^*E=\cO_N\otimes_{\cO_M}E$,
\item
for any section $\eta$ of $\Theta_N$, $(f^*\nabla)_\eta=(\cL_\eta\otimes\id)+\nabla_{Tf(\eta)}$ and $(f^*\Phi)_\eta=\Phi_{Tf(\eta)}$,
\item
$f^*R_\infty=\id\otimes R_\infty$, $f^*R_0=\id\otimes R_0$,
\end{itemize}
where $\nabla$ and $\Phi$ are understood to be linearly extended to $f^*\Theta_M$, and $\cL_\eta$ denotes the Lie derivative with respect to $\eta$.

If $(M,E,\nabla,R_\infty,\Phi,R_0)$ is pre-Saito structure on $M$, then so is its pull-back on $N$. If $f$ is a closed immersion, then we say that $(\bE,\nablab)$ is a \emph{deformation} of $f^*(\bE,\nablab)$.

\begin{example}[of a deformation]\label{exam:deformation}
Let us start with a pre-Saito structure on a point, that is, a triple $(E^o,R_\infty,R_0^o)$, where $E^o$ is a finite dimensional vector space and $R_\infty,R_0^o$ are two endomorphisms of $E^o$. We consider the following ``trivial'' one-parameter deformation $(\Afu,E=\cO_{\Afu}\otimes_\CC E^o,\nabla,R_\infty,\Phi,R_0)$ (parametrized by the complex line $\Afu$ with coordinate $x$), with:
\begin{align*}
\nabla&=d,\\
R_\infty&=\id\otimes R_\infty,\\
R_0(x)&=e^{x(R_\infty+\id)}\cdot (\id\otimes R_0^o) \cdot e^{-xR_\infty}=e^{x(\ad R_\infty+\id)}(\id\otimes R_0^o),\\
\Phi&=-R_0(x)dx.
\end{align*}
The only non-trivial relation to be checked~is
\[
\Phi_{\partial_x}+\dfrac{\partial R_0}{\partial x}+[R_\infty,\Phi_{\partial_x}]=0,
\]
which follows from the definition of $R_0$, as $\Phi_{\partial_x}=-R_0$. Let us remark that, according to this relation, any one-parameter deformation with $\Phi_{\partial_x}=-R_0(x)$ is isomorphic to the previous one.

One can also remark that the eigenvalues of $R_0(x)$ are $e^x$ times the eigenvalues of $R_0^o$.

Last, let us notice that, if $R_\infty$ is semisimple with integral eigenvalues, we can define the family in an algebraic way with respect to the variable $\lambda\in\CC^*$, by replacing $e^x$ with~$\lambda$.
\end{example}

\begin{remarks}[on Example \ref{exam:deformation}]\label{rem:deformation}\mbox{}
\begin{enumerate}
\item\label{rem:deformation1}
From the point of view of the data $(\bE,\nablab)$, the construction of Example \ref{exam:deformation} only consists of a rescaling in the variable $\hb$. On $\bE^o$ we have the connection $\nablab^o=d_\hb+(R_\infty-\hb R_0^o)d\hb/\hb$ and, if $E=\CC[\lambda,\lambda^{-1}]\otimes_\CC E^o$, we consider on $\bE$ the trivial connection $\nablab'=d_\lambda+\nablab^o=d+(R_\infty-\hb R_0^o)d\hb/\hb$. Let us now consider the rescaling
\[
\rho^*:\CC[\lambda,\lambda^{-1},\hb]\to\CC[\lambda,\lambda^{-1},\hb],\qquad \lambda\mto \lambda,\quad \hb\mto \lambda\hb.
\]
The inverse image of $\nablab'$ by this rescaling is $d+(R_\infty-\lambda\hb R_0^o)(d\lambda/\lambda+d\hb/\hb)$. It has Poincar\'e rank one along $\hb=\infty$ (we are not interested in the behaviour when $\lambda\to0$ or $\lambda\to\infty$). Up to now, the construction is algebraic. However, we need to change the trivialization so that $\nablab$ gets the Birkhoff normal form. In order to do so, we pull-back $(\bE,\nablab)$ by the uniformization $\CC\to\CC^*$, $x\mto\lambda=e^x$, and we change the trivialization using $e^{xR_\infty}$. Let us also notice that the uniformization $\lambda=e^x$ is not needed if $R_\infty$ is semisimple with integral eigenvalues.
\item\label{rem:deformation2}
The construction of Example \ref{exam:deformation} can be done starting from any pre-Saito structure $(M,\cO_M\otimes_\CC E^o,d,R_\infty,\Phi,R_0)$ to produce a pre-Saito structure $(M\times\nobreak\Afu,\cO_{M\times\Afu}\otimes_\CC E^o,d,R_\infty,\wt\Phi,\wt R_0)$ with
\begin{align*}
\wt R_0&=e^{x(\id+\ad R_\infty)}(\id\otimes R_0),\\
\wt\Phi&=e^{x(\id+\ad R_\infty)}(\id\otimes\Phi)-\wt R_0 dx.
\end{align*}
If the kernel of $\id+\ad R_\infty$ is non zero, then there could exist other pre-Saito structures (\ie other $\wt\Phi$) with the same $\wt R_0$.
\item\label{rem:deformation3}
The construction of Example \ref{exam:deformation} can be iterated, using \eqref{rem:deformation2}, but this does not lead to any interesting new deformation.
\end{enumerate}
\end{remarks}

Let $i:M\to N$ be an immersion. We say that a pre-Saito structure $(\bE_N,\nablab_N)$ on $N$ is a \emph{universal deformation} of its restriction $(\bE,\nablab)\defin i^*(\bE_N,\nablab_N)$ if any other deformation of $(\bE,\nablab)$ comes from $(\bE_N,\nablab_N)$ by a unique base change inducing the identity on $M$.

If $(M,E,\nabla)$ is a vector bundle with flat connection, then there is no loss of information by fixing a horizontal trivialization $(E,\nabla)\simeq(\cO_M\otimes_\CC E^o,d)$, where $E^o=\ker\nabla$ is the space of $\nabla$-horizontal sections of~$E$, which can also be identified with $E/\gm E$.

Let $(M,\cO_M\otimes_\CC E^o,d,R_\infty,\Phi,R_0)$ be a pre-Saito structure. Let $\omega^o$ be any element of $E^o$ and let $\omega=1_M\otimes\omega^o$ denote the unique $\nabla$-horizontal section determined by~$\omega^o$. Then $\Phi$ defines a morphism $\varphi_\omega:\Theta_M\to \cO_M\otimes_\CC E^o$, $\xi\mto-\Phi_\xi(\omega)$, which can be regarded as a section of $\Omega^1_M\otimes_\CC E^o$, and which is called the \emph{infinitesimal period mapping attached to $\omega^o$}. The conditions $d\Phi=0$ and $d\omega=0$ imply $d\varphi_\omega=0\in\Omega^2_M\otimes_\CC E^o$.

\begin{proposition}[Hertling-Manin \cite{H-M04}]\label{prop:herman}
Let $(M,\cO_M\otimes_\CC E^o,d,R_\infty,\Phi,R_0)$ be a germ at~$o\in M$ of pre-Saito structure. Let us assume that there exists $\omega^o\in E^o$ such that $\omega^o$ and its images under the iteration of the maps $R_0^o:E^o\to E^o$ and $\Phi^o_\xi:E^o\to E^o$ (for all $\xi\in\Theta_M^o$) generate $E^o$. Let us set $\omega=1_M\otimes\omega^o$.

Let $N$ be a germ of complex analytic manifold along $M$ and let $i:M\hto N$ denote the immersion. Then, there is a one-to-one correspondence between deformations $(N,\cO_N\otimes_\CC E^o,d,R'_\infty,\Phi',R'_0)$ of the pre-Saito structure $(M,\cO_M\otimes_\CC E^o,d,R_\infty,\Phi,R_0)$ parametrized by $N$ and germs $\varphi\in\Omega^1_N\otimes_\CC E^o$ such that
\[\tag{$*$}
\begin{cases}
i^*\varphi=\varphi_\omega,\\
d\varphi=0,
\end{cases}
\]
the correspondence being given by
\[\tag{$**$}
(N,\cO_N\otimes_\CC E^o,d,R'_\infty,\Phi',R'_0)\mto\varphi=\varphi_{1_N\otimes\omega^o}.
\]
\end{proposition}

\begin{proof}
We set $\cO_M=\CC\{x\}$ with $x=(x_1,\dots,x_m)$ and $\cO_N=\CC\{x,y\}$ with $y=(y_1,\dots,y_n)$. On the one hand, it is easy to check that $\varphi$ defined by $(**)$ satisfies Properties $(*)$. Let us thus start, on the other hand, with $\varphi$ satisfying $(*)$. Clearly, if a deformation exists, then $1_N\otimes\omega^o$ is horizontal. We can therefore argue by induction on $n$ and assume that $n=1$. We will thus set $y=y_1$. Let us also remark that, under the assumption on $\omega^o$, the images of $\omega$ under the iteration of the maps $R_0:\cO_M\otimes_\CC E^o\to \cO_M\otimes_\CC E^o$ and $\Phi_\xi:\cO_M\otimes_\CC E^o\to \cO_M\otimes_\CC E^o$ generate $\cO_M\otimes_\CC E^o$ as a $\cO_M$-module.

Let us fix a basis $\bme^o$ of $E^o$. We then get matrices $C^{(i)}(x)$, $B_0(x)$ and $B_\infty$ satisfying \eqref{eq:relprefrobhoriz}. If the desired pre-Saito structure exists, it must have $1\otimes \bme^o$ as horizontal basis. So, we search for $C^{\prime(i)}(x,y)$, $D'(x,y)$, $B'_0(x,y)$ (and we set $B'_\infty=B_\infty$) satisfying \eqref{eq:relprefrobhoriz} with one variable more (where $D'$ is the component of $C'$ on $dy$). One sets $C^{\prime(i)}(x,y)=\sum_{k\geq0}C^{\prime(i)}_k(x)y^k$, etc.\ and one computes inductively the coefficients $C^{\prime(i)}_k(x)$, $D'_k(x)$, $B'_{0,k}(x)$.

One sets first $C^{\prime(i)}_0(x)=C^{(i)}(x)$ and $B'_{0,0}(x)=B_0(x)$. One also must have $\sum_iC^{\prime(i)}(\omega^o)dx_i+D'(\omega^o)dy=-\varphi$.

If $C^{\prime(i)}_{\leq k}(x)$, $B'_{0,\leq k}(x)$ and $D'_{\leq k-1}(x)$ are found (satisfying \eqref{eq:relprefrobhoriz} mod~$y^k$), the generating assumption and the desired commutation of $D'$ with $C^{\prime(i)}$ and $B'_0$ implies that $D'_{\leq k}(\omega^o)$ (which is determined by $\varphi$, hence known) uniquely determines such a~$D'_{\leq k}$. Let us also notice for future use that, modulo $y^{k+1}$, $D'_{\leq k}$ belongs then to the commutative algebra generated by the classes modulo $y^{k+1}$ of the $C^{\prime(i)}_{\leq k}$ and $B'_{0,\leq k}(x)$.

Then $C^{\prime(i)}_{\leq k+1}$ and $B'_{0,\leq k+1}$ are uniquely determined by their initial value and the equations
\[
\frac{\partial C^{\prime(i)}_{\leq k+1}}{\partial y}=\frac{\partial D'_{\leq k}}{\partial x_i},\qquad
\frac{\partial B'_{0,\leq k+1}}{\partial y}= [B_\infty,D'_{\leq k}]-D'_{\leq k}.
\]
That all desired relations at the level $k+1$ are satisfied is then easily verified. It remains to prove convergence. This is done in \cite{H-M04}.
\end{proof}

\begin{remark}\label{rem:herman}
Let us assume that the conditions of the proposition are fulfilled. Given $\varphi_\omega\in\Omega^1_M\otimes_\CC E^o$ satisfying $d\varphi_\omega=0$, an extension $\varphi=\sum_i\varphi_idx_i+\sum_j\psi_jdy_j$ as in the proposition is determined in a unique way from $\psi=\sum_j\psi_jdy_j$ provided that $\psi$ is $d_y$-closed. Therefore, there is a one-to-one correspondence between the deformations $(N,\cO_N\otimes_\CC E^o,d,R_\infty,\Phi',R'_0)$ as in the proposition (with a chosen projection $N\to M$), and the set of $\Psi\in\cO_N\otimes_\CC E^o$ satisfying $\Psi(x,0)=0$: one associates to $\Psi$ the unique $\Phi',R'_0$ defined by $\varphi$, where $\varphi$ is determined by $\psi=d_y\Psi$.

In particular, if we fix $\wt\varphi=\sum_i\varphi_{\omega,i}(x)dx_i+\sum_j\psi_j(x,0)dy_j$, that is, if we fix $\psi(x,0)$, there exists $\psi(x,y)$ which is $d_y$-closed and restricts to $\psi(x,0)$ at $y=0$. Therefore, given any such $\wt\varphi$, there exists a deformation $(N,\cO_N\otimes_\CC E^o,d,R_\infty,\Phi',R'_0)$ such that $\varphi_{1\otimes\omega^o\mid y=0}=\wt\varphi$.
\end{remark}

\begin{corollary}[Hertling-Manin \cite{H-M04}]\label{cor:hermanuniv}
Let $(M,E,\nabla,R_\infty,\Phi,R_0)$ be a germ of pre-Saito structure with $\omega^o\in E^o$ satisfying the assumptions of Proposition \ref{prop:herman}. If moreover $\varphi_\omega^o:\Theta_M^o\to E^o$ is injective, then $(M,E,\nabla,R_\infty,\Phi,R_0)$ has a universal deformation parametrized by the germ $\wt E^o=(E^o,0)$.
\end{corollary}

\begin{proof}
As above, let us identify $(E,\nabla)$ with $(\cO_M\otimes_\CC E^o,d)$. For any $N$ as above, we can therefore identify a section of $\cO_N\otimes_\CC E^o$ vanishing at $o$ with a morphism $N\to \wt E^o$, where $\wt E^o$ is the analytic germ of the $\CC$-vector space $E^o$ at the origin. For $\varphi$ as in the proposition, we have $d\varphi=0$, hence $\varphi=d\chi$ where $\chi\in\cO_N\otimes_\CC E^o$ is uniquely determined by the initial condition $\chi(o)=0$. We regard $\chi$ as a morphism $\chi:N\to \wt E^o$. In particular, to $\varphi_\omega$ we associate $\chi_\omega:M\to \wt E^o$.

From Proposition \ref{prop:herman}, one deduces that giving a deformation of the pre-Saito structure $(M,E,\nabla,R_\infty,\Phi,R_0)$ parametrized by $N$ is equivalent to giving a commutative diagram
\[
\xymatrix@C=1.5cm{
M\ar[r]^-{\chi_\omega}\ar@{_{ (}->}[d]&\wt E^o\\
N\ar[ru]^-{\chi}&
}
\]
and, given a base change $\nu:N'\to N$ inducing the identity on $M$, the pull-back $\nu^*(M,E,\nabla,R_\infty,\Phi,R_0)$ corresponds to $\chi'=\chi\circ\nu$. In particular, $(N,\chi)$ is universal if and only if for any $(N',\chi')$ there exists a unique $\nu:N'\to N$ inducing the identity on $M$ such that $\chi'=\chi\circ\nu$. The assumption on $\varphi_\omega^o$ means that $\chi_\omega$ is an immersion. The universal deformation must then correspond to the diagram
\[
\xymatrix@C=1.5cm{
M\ar@{^{ (}->}[r]^-{\chi_\omega}\ar@{_{ (}->}[d]_{\chi_\omega}&\wt E^o\\
\wt E^o\ar[ru]^-{\id}
}
\]
\end{proof}

From the last point in Remark \ref{rem:herman}, we get:

\begin{corollary}\label{cor:univphio}
Let $(M,E,\nabla,R_\infty,\Phi,R_0)$ be a germ of pre-Saito structure with $\omega^o\in E^o$ satisfying the assumptions of Proposition \ref{prop:herman}. Given any smooth analytic germ $N\supset M$ together with an isomorphism $\wt\varphi:i^*\Theta_N\to E$ restricting to $\varphi_\omega$ on $\Theta_M\subset i^*\Theta_N$, there exists on $N$ a universal deformation of $(M,E,\nabla,R_\infty,\Phi,R_0)$ such that $\varphi_{1\otimes\omega\mid M}=\wt\varphi$.\qed
\end{corollary}

Such a deformation is not unique, but one can obtain any such deformation from a given one through a unique base change $N\to N$, and it is tangent to the identity when restricted to $M$.

Concerning uniqueness, one also obtains:

\begin{corollary}\label{cor:hermanuniqueness}
Under the assumptions of Corollary \ref{cor:hermanuniv}, let us consider two deformations of $(M,E,\nabla,R_\infty,\Phi,R_0)$ with parameter spaces $N,N'\supset M$ being two smooth analytic germs, for which the corresponding $\chi,\chi':N,N'\to \wt E^o$ are immersions with the same image. Then these two deformations are isomorphic, \ie one comes from the other by a base change inducing an isomorphism on tangent bundles.\qed
\end{corollary}

\begin{example}\label{exam:deformationgenerating}
In the situation of Example \ref{exam:deformation}, let us assume that $R_0^o$ has a cyclic vector $\omega^o\in E^o$. If $d=\dim_\CC E^o$, then $\omega^o,\dots,(R_0^o)^{d-1}(\omega^o)$ is a basis of~$E^o$. The generating condition of Proposition \ref{prop:herman} is satisfied, hence there exists a universal deformation of the one-parameter pre-Saito structure defined in Example \ref{exam:deformation}, parametrized by $\wt E^o$. Setting $\omega=1\otimes\omega^o$, and denoting by $x_1$ the coordinate on $\Afu$, the map
\[
\varphi_\omega:\Theta_{\Afu}\to\cO_{\Afu}\otimes_\CC E^o
\]
is given by
\[
\varphi_\omega(\partial_{x_1})=-\Phi_{\partial_{x_1}}(\omega)=R_0(x_1)(\omega^o)=e^{x_1(\id+\ad R_\infty)}(R_0^o)(\omega^o),
\]
and we have
\[
\chi_\omega(x_1)=\frac{e^{x_1(\id+\ad R_\infty)}-\id}{(\id+\ad R_\infty)}\,(R_0^o)(\omega^o).
\]
Using the notation of Remark \ref{rem:herman}, let us set $\wt\varphi=\sum_{j=0}^{d-1}R_0(x_1)^j(\omega^o)dx_j$. It defines an isomorphism $\Theta_{\Afu\times(\CC^{d-1},0)|\Afu}\to\cO_{\Afu}\otimes_\CC E^o$, and induces a local biholomorphic map $\chi:\Afu\times(\CC^{d-1},0)\to E^o$. If $(x_0,\dots,x_{d-1})$ denote the coordinates on $\Afu\times(\CC^{d-1},0)$, we thus have $\Phi_{\partial_{x_j}}(1\otimes\omega^o)_{|\Afu}=R_0(x_1)^j(\omega^o)$.

Let us consider the case where $\omega^o$ and $R_0^o(\omega^o)$ are eigenvectors of $R_\infty$ with respective eigenvalues $\delta_0,\delta_1$. Then
\[
\varphi_\omega(\partial_{x_1})=e^{(\delta_1-\delta_0+1)x_1}R_0^o(\omega^o)\quad\text{and}\quad \chi_\omega(x_1)=\frac{e^{(\delta_1-\delta_0+1)x_1}-1}{\delta_1-\delta_0+1}\, R_0^o(\omega^o).
\]
In such a case $\chi_\omega$ is a parametrization of the line $\CC\cdot R_0^o(\omega^o)$ minus the point $(-1/(\delta_1-\delta_0+1))R_0^o(\omega^o)$ (if $\delta_1=\delta_0-1$, this point is at infinity, so does not have to be deleted). Moreover, for any $x^o\in\Afu$, the analytic germ $(E^o,\chi_\omega(x^o))$ is the universal deformation of the germ at $x^o$ of the pre-Saito structure constructed in Example \ref{exam:deformation}. In the local coordinates $(x_0,\dots,x_{d-1})$, we have $\Phi_{\partial_{x_j}|\Afu}=R_0(x_1)^j$: indeed, this holds when applying both operators to $\omega^o$; as $\omega^o$ is a cyclic vector for $R_0(x_1)$ for any $x_1$, and as $\Phi_{\partial_{x_j}|\Afu}$ commutes with $R_0(x_1)$, we get the desired assertion.

If $\delta_1=\delta_0-1$, then $\chi_\omega(x_1)=x_1R_0^o(\omega^o)$ defines a closed embedding $\Afu\hto E^o$. We will mainly consider this case later on, and we will then denote by $\wt E^o$ the analytic germ $(E^o,\CC\cdot R_0^o(\omega^o))$.
\end{example}

\subsection{Pre-Saito structures with a finite group action}\label{subsec:W}
Let $M$ be a punctual germ of complex manifold and let us assume that $M$ is acted on by a finite group $\rW$ of automorphisms. For $\rw\in\rW$, we denote by $\rw:M\to M$ the corresponding automorphism and by $\rw^*:\cO_M\to\cO_M$ the associated morphism of $\CC$-algebras. The fixed subspace $M^W$ is also a smooth analytic germ.

If $E$ is a free $\cO_M$-module, we say that the action of $\rW$ lifts linearly to $E$ if, for any $\rw\in\rW$, there exists an isomorphism $a_\rw:E\to\rw^*E$ and, for any $\rw,\rw'\in \rW$, the following diagram commutes:
\[
\xymatrix{
E\ar[rd]_{a_\rw}\ar[rr]^-{a_{\rw'\rw}}&&\rw^*\rw'^*E\\
&\rw^*E\ar[ur]_-{\rw^*a_{\rw'}}&
}
\]
In particular, the restriction $E_{|M^W}$ is equipped with a linear action of $W$. For instance, there is a canonical linear lifting of the $\rW$-action to the tangent bundle $\Theta_M$, if we set $a_\rw=T\rw:\Theta_M\to\rw^*\Theta_M$, and we have $\Theta_{M^W}=(\Theta_{M|M^W})^W$.

Given a pre-Saito structure $(M,E,\nabla,R_\infty,\Phi,R_0)$, we say that the $\rW$-action on~$M$ lifts linearly to $(M,E,\nabla,R_\infty,\Phi,R_0)$ if it lifts linearly to $E$ and each $a_\rw$ induces an isomorphism of pre-Saito structures. If we fix the horizontal trivialization $(E,\nabla)\simeq(\cO_M\otimes_\CC E^o,d)$, then we must have $a_\rw=\id\otimes a^o_\rw$, with $a_{\rw'\rw}^o=a_{\rw'}^oa_\rw^o$ for any $\rw,\rw'\in \rW$. If we fix coordinates $(x_1,\dots,x_m)$ on $M$, we get, setting $\rw_j=x_j\circ\rw$,
\begin{align*}
a^o_\rw R_\infty(a^o_\rw)^{-1}&=R_\infty,\\
a^o_\rw R_0(x)(a^o_\rw)^{-1}&=R_0(\rw(x)),\\
a^o_\rw \Phi_{\partial_{x_i}}(x)(a^o_\rw)^{-1} &=\sum_k\frac{\partial\rw_k}{\partial x_i}(x)\cdot \Phi_{\partial_{x_k}}(\rw(x)).
\end{align*}
In particular, $\rW$ acts $\CC$-linearly on $E^o$ and $a_\rw^o$ commutes with $R_\infty$ and $R_0^o$.

Let $\omega=\sum_i\omega_i(x)\otimes e_i^o$ be a section of $E$. We have $\rw^*\omega=\sum_i\omega_i(\rw(x))\otimes e_i^o$ and $a_\rw(\omega)=\sum_i\omega_i(x)\otimes a_\rw^o(e_i^o)$. We say that $\omega$ is $\rW$-equivariant if, for any $\rw\in \rW$, we have $\rw^*(\omega)=a_\rw(\omega)$. If $\omega$ is $\rW$-equivariant, then its restriction to $M^\rW$ is $\rW$-invariant. Conversely, assume that $\omega_{M^\rW}$ is a flat $\rW$-invariant section of $E_{|M^\rW}$. Let $\omega$ be its flat extension to $E$. Then $\omega$ is $\rW$-equivariant. Similarly, if $\omega^o\in E^o$ is $\rW$-invariant, then its flat extension $\omega$ is $\rW$-equivariant.

If $\omega$ is $\rW$-equivariant, then the following diagram commutes:
\begin{equation}\label{eq:phiequiv}
\begin{array}{c}
\xymatrix@C=1.5cm{
\Theta_M\ar[r]^-{\varphi_\omega}\ar[d]_{T\rw}&E\ar[d]^{a_\rw}\\
\rw^*\Theta_M\ar[r]^-{\rw^*\varphi_\omega}&\rw^*E
}
\end{array}
\end{equation}
Moreover, $\wt E^o$ is naturally equipped with a $\rW$-action (coming from the linear action on $E^o$) and $\chi_\omega:M\to\wt E^o$ is $\rW$-equivariant.

\subsubsection*{$\rW$-equivariant version of Proposition \ref{prop:herman} and Corollary \ref{cor:hermanuniv}}
Let $(M,E,\nabla,R_\infty,\Phi,R_0)$ be a pre-Saito structure with $\rW$-action and let $\omega^o\in E^o$. Let us assume that $\omega^o$ is $\rW$-invariant and let $\omega$ be its flat extension, which is $\rW$-equivariant.

By a $\rW$-equivariant deformation of $(M,E,\nabla,R_\infty,\Phi,R_0)$ we mean a deformation parametrized by $N\supset M$ with a $\rW$-action, such that $M$ is left stable by the $\rW$-action on $N$, and which restricts (with $\rW$-action) to $(M,E,\nabla,R_\infty,\Phi,R_0)$. Proposition \ref{prop:herman} extends as follows:

\begin{corollary}\label{cor:hermanWequiv}
With the assumptions of Proposition \ref{prop:herman}, let us moreover assume that
\begin{enumerate}
\item
$\omega^o$ is $\rW$-invariant,
\item
the $\rW$-action on $M$ extends to a $\rW$-action on $N$.
\end{enumerate}
Then, under the correspondence of Proposition \ref{prop:herman}, $\rW$-equivariant deformations of $(M,E,\nabla,R_\infty,\Phi,R_0)$ parametrized by $N$ correspond to $\rW$-equivariant closed sections $\varphi\in\Omega_N^1\otimes_\CC E^o$ (\ie the diagram corresponding to \eqref{eq:phiequiv} commutes).
\end{corollary}

\begin{proof}
Starting from a $\rW$-equivariant deformation, and as $\omega^o$ is $\rW$-invariant, the associated $\varphi$ is easily seen to be $\rW$-equivariant. Conversely, if $\varphi$ is $\rW$-equivariant, then we have two deformations of $(M,E,\nabla,R_\infty,\Phi,R_0)$ defined in coordinates $x$ on $N$ by $R_\infty$, $R_0(x)$  and $\Phi_{\partial_{x_i}}(x)$ on the one hand, and by $(a_\rw^o)^{-1}R_\infty a_\rw^o$, $(a_\rw^o)^{-1}R_0(\rw(x)) a_\rw^o$ and $\sum_k\frac{\partial\rw_k}{\partial x_i}(x)\cdot (a_\rw^o)^{-1}\Phi_{\partial_{x_k}}(\rw(x))a_\rw^o$ on the other hand. That $\varphi$ is $\rW$-equivariant means that the $\varphi_\omega$ associated to each of these deformations coincide. By uniqueness in Proposition \ref{prop:herman}, these deformations coincide.
\end{proof}

\begin{corollary}\label{cor:hermanWuniv}
With the assumptions of Corollary \ref{cor:hermanWequiv}, let us moreover assume that $\varphi_\omega^o:\Theta_M^o\to E^o$ is an immersion. Then the $\rW$-action on $\wt E^o$ coming from the linear action on $E^o$ can be lifted as a $\rW$-action on the universal deformation of $(M,E,\nabla,R_\infty,\Phi,R_0)$ given by Corollary \ref{cor:hermanuniv} and this action restricts, through $\chi_\omega$, to the given one on $(M,E,\nabla,R_\infty,\Phi,R_0)$.\qed
\end{corollary}

\subsubsection*{$\rW$-invariant version of Proposition \ref{prop:herman} and Corollary \ref{cor:hermanuniv}}
Let us consider the case where $\rW$ acts trivially on $M$. We will say that the action is linearized. In particular, $E^o$ has a $\rW$-action and $R_0^o$ and the $\Phi_\xi^o$ commute with this $\rW$-action.

One can define the notion of a deformation with linearized $\rW$-action, and that of a universal deformation with linearized $\rW$-action. The results of Hertling and Manin extend as follows:
\begin{enumerate}
\item
In Proposition \ref{prop:herman}, one assumes that $\omega^o$ is $\rW$-invariant and that the images of $\omega^o$ under the iteration of the maps $R_0^o$ and $\Phi_\xi^o$ generate the invariant subspace $(E^{o})^\rW$. Then $\omega$ is a section of $E^\rW$. The $\rW$-invariant version of Proposition \ref{prop:herman} is that there exists a one-to-one correspondence between the set of deformations with linearized $\rW$-action and the set of $\varphi\in\Omega^1_N\otimes_\CC (E^o)^\rW$ satisfying $(*)$.

For the proof, one notices that, by induction, the matrices $C_k^{\prime(i)}$, $D'_k$ and $B'_{0,k}$ commute with the $\rW$-action.
\item
If moreover $\varphi_\omega^o$ is an immersion $\Theta_M^o\hto(\wt E^o)^\rW$, then $(\wt E^o)^\rW$ is the base space of a universal deformation with linearized $\rW$-action of $(M,E,\nabla,R_\infty,\Phi,R_0)$ (same proof as for Corollary \ref{cor:hermanuniv}).
\item
Corollaries \ref{cor:univphio} and \ref{cor:hermanuniqueness} extend in the same way.
\end{enumerate}

\begin{remark}
Let $(M,E,\nabla,R_\infty,\Phi,R_0)$ be a pre-Saito structure with a (not linearized) $\rW$-action, and let $\omega^o\in E^o$ be $\rW$-invariant. Let us assume that $\omega^o$ fulfills the conditions in Corollary \ref{cor:hermanuniv}. Then the universal deformation with parameter space $\wt E^o$ comes equipped with a (non linearized) $\rW$-action. The restriction of this deformation to the subspace $(\wt E^o)^\rW$ has therefore a linearized $\rW$-action. However, it may not be, as such, a universal deformation with linearized $\rW$-action of $(M,E,\nabla,R_\infty,\Phi,R_0)_{|M^\rW}$, as the images of $\omega^o$ under the iterates of $R_0^o$ and the $\Phi^o_\xi$ ($\xi\in\Theta_{M^\rW}^o$) may not generate $(E^o)^\rW$. One can ask whether there exists an interesting smooth subspace contained in $(\wt E^o)^\rW$ so that $(M,E,\nabla,R_\infty,\Phi,R_0)_{|M^\rW}$ is the universal deformation with linearized action of its restriction to this subspace.
\end{remark}

\subsection{Metric}

\begin{definition}\label{def:prefrobwith}
A pre-Saito structure $(M,E,\nabla,R_\infty,\Phi,R_0,g)$ (with metric) of weight $w$ consists of the following data:
\begin{enumerate}
\item
A pre-Saito structure (without metric) $(M,E,\nabla,R_\infty,\Phi,R_0)$ as in \S\ref{subsec:prefro},
\item
a nondegenerate symmetric $\cO_M$-bilinear form $g$ on $E$
\end{enumerate}
which satisfy the following relations, denoting by ${}^*$ the adjoint with respect to $g$:
\[
\nabla(g)=0,\quad R_\infty+R_\infty^*=-w\id,\quad \Phi^*=\Phi,\quad R_0^*=R_0.
\]
\end{definition}

Let us notice that $\Phi^*=\Phi$ means that for all $\xi\in\Theta_M$, $(\Phi_\xi)^*=\Phi_\xi$.

\begin{example}\label{exam:deformationmetric}
In the situation of Example \ref{exam:deformation}, if we moreover have a nondegenerate symmetric bilinear form $g^o$ on $E^o$ such that $R_0^{o*}=R_0^o$ and $R_\infty^*+R_\infty=-w\id$, then, setting $g=\id\otimes g^o$ (\ie extending trivially $g^o$, so that $g$ is $\nabla$-flat) we still have $R_\infty^*+R_\infty=-w\id$ and
\[
R_0^*=e^{-xR_\infty^*}R_0^{o*}e^{(xR_\infty^*+\id)}=e^{xR_\infty}\,e^{wx\id}R_0^o\,e^{(1-w)x\id}e^{-xR_\infty}=R_0,
\]
so the deformed pre-Saito structure remains of weight $w$.
\end{example}

\begin{corollary}[Hertling-Manin \cite{H-M04}]\label{cor:hermanunivwith}
Let $(M,\cO_M\otimes_\CC E^o,d,R_\infty,\Phi,R_0)$ be a pre-Saito structure with $\omega^o\in E^o$ satisfying the conditions of Proposition \ref{prop:herman}. Let $(N,\cO_N\otimes_\CC E^o,d,R_\infty,\Phi',R'_0)$ be any deformation of $(M,\cO_M\otimes_\CC E^o,d,R_\infty,\Phi,R_0)$. Assume that $g$ is a flat metric on $\cO_M\otimes_\CC E^o$ giving $(M,\cO_M\otimes_\CC E^o,d,R_\infty,\Phi,R_0)$ weight $w$ and let $g'$ be the unique $d$-flat metric on $\cO_N\otimes_\CC E^o$ extending~$g$.

Then $(N,\cO_N\otimes_\CC E^o,d,R_\infty,\Phi',R'_0,g')$ is a pre-Saito structure of weight $w$.
\end{corollary}

\begin{proof}
In the proof of Proposition \ref{prop:herman}, let us choose the basis $\bme^o$ so that it is orthonormal with respect to $g^o$. Assume, by induction, that the matrices $C_{\leq k}^{\prime(i)}$, $B'_{0,\leq k}$ and $D'_{\leq k-1}$ are symmetric. Then $D'_{\leq k}$ is symmetric, as it can be expressed as a polynomial in $C_{\leq k}^{\prime(i)},B'_{0,\leq k}$ modulo $y^{k+1}$, then $\partial C_{\leq k+1}^{\prime(i)}/\partial y$ and $\partial B'_{0,\leq k+1}/\partial y$ are symmetric, hence also $C_{\leq k+1}^{\prime(i)},B'_{0,\leq k+1}$.
\end{proof}

\begin{remark}
The adaptation of the previous result with $\rW$-action is straightforward.
\end{remark}

\subsection{Frobenius manifolds}

We still assume that $M$ is a punctual analytic germ. Let us recall well-known results (see \eg \cite[Chap.~VII]{Bibi00}).
\begin{definition}\label{def:Frobenius}
A Frobenius manifold structure $(M,\star,g,e,\gE)$ of weight~$D$ consists of
\begin{enumeratei}
\item
A symmetric nondegenerate $\cO_M$-bilinear form $g$ on $\Theta_M$, with associated Levi-Civita (\ie torsion free) connection $\nabla:\Theta_M\to\Omega_M^1\otimes_{\cO_M}\Theta_M$;
\item
A $\cO_M$-bilinear product $\star$ on $\Theta_M$;
\item
Two sections $e$ and $\gE$ of $\Theta_M$;
\end{enumeratei}
\noindent
subject to the following relations:
\begin{enumeratea}
\item
$\nabla$ is \emph{flat};
\item
$\star$ is commutative and associative;
\item
$e$ is a unit for $\star$ and is $\nabla$-horizontal;
\item
$\cL_{\gE}(e)=-e$, $\cL_{\gE}(\star)=\star$, $\cL_{\gE}(g)=Dg$ for some $D\in\CC$;
\item
If $c\in\Gamma(M,\Omega_M^1{}^{\otimes 3})$ is defined by $c(\xi_1,\xi_2,\xi_3)=g(\xi_1\star\xi_2,\xi_3)$, then $\nabla c$ is symmetric in its four arguments.
\end{enumeratea}
\end{definition}

Let $(M,E,\nabla,R_\infty,\Phi,R_0,g)$ be a pre-Saito structure. Let $\omega$ be a $\nabla$-horizontal section of $E$. It defines a $\cO_M$-linear morphism $\varphi_\omega:\Theta_M\to E$ by $\xi\mto-\Phi_\xi(\omega)$.

\begin{definition}
Given a pre-Saito structure $(M,E,\nabla,R_\infty,\Phi,R_0,g)$ of weight~$w$, we say that a $\nabla$-horizontal section $\omega$ of $E$ is
\begin{enumerate}
\item
\emph{primitive} if the associated period mapping $\varphi_\omega:\Theta_M\to E$ is an isomorphism,
\item
\emph{homogeneous} of degree $q\in\CC$ if $R_\infty\omega=q\omega$.
\end{enumerate}
A pre-Saito structure $(M,E,\nabla,R_\infty,\Phi,R_0,g)$ of weight~$w$ equipped with a primitive homogeneous section $\omega$ is called a \emph{Saito structure}.
\end{definition}

If $\omega$ is primitive and homogeneous, $\varphi_\omega$ induces a flat torsion free connection $\onabla$ on $\Theta_M$, and an associative and commutative $\cO_M$-bilinear product $\star$, with $e=\varphi_\omega^{-1}(\omega)$ as unit, and $\onabla e=0$.

The Euler field is $\gE=\varphi_{\omega}^{-1}(R_0(\omega))$. It is therefore a section of $\Theta_M$. We have $\onabla\gE={}^\omega\!R_\infty+q\id$, with ${}^\omega\!R_\infty=\varphi_\omega^{-1}\circ R_\infty\circ\varphi_\omega-\id$. In particular, $\onabla\onabla\gE=0$.

\begin{remark}
We have $\cL_{\gE}(e)=-e$, $\cL_{\gE}(\star)=\star$. If we set $D=2q+2-w$, we have, for ${}^\omega\!g$ induced by $\varphi_\omega$ as above:
\[
\cL_{\gE}({}^\omega\!g)= D\cdot {}^\omega\!g.
\]
\end{remark}

\begin{proposition}
Let $(M,E,\nabla,R_\infty,\Phi,R_0,g)$ be a pre-Saito structure of weight~$w$. To any homogenous primitive section $\omega$ of $E$ having weight $q$ is associated canonically on $M$, through the infinitesimal period mapping $\varphi_\omega$, a Frobenius manifold structure of weight $D=2q+2-w$.

Conversely, any Frobenius manifold structure $(M,\star,g,e,\gE)$ defines a Saito structure $(M,E,\nabla,R_\infty,\Phi,R_0,g)$ having $e$ as homogeneous primitive form.
\end{proposition}

\begin{proof}
Let us give the correspondence $(M,\star,g,e,\gE)\mto(M,E,\nabla,R_\infty,\Phi,R_0,g)$. We define:
\begin{itemize}
\item
$E=\Theta_M$,
\item
$\nabla$ is the Levi-Civita connection of $g$,
\item
$\Phi_\xi(\eta)=-(\xi\star\eta)$,
\item
$R_0=\gE\star{}=-\Phi_\gE$,
\item
$R_\infty=\nabla\gE-\id$,
\item
$q=0$, $w=2-D$.\qedhere
\end{itemize}
\end{proof}

\begin{remark}\label{rem:sqrt}
Let $(M,E,\nabla,R_\infty,\Phi,R_0,g)$ be a pre-Saito structure. If $\omega$ is a primitive homogeneous section of $E$, then so is $\lambda\omega$ for any $\lambda\in\CC^*$. It gives rise to the Frobenius manifold structure $(M,\star,\lambda^2g,e,\gE)$. In particular, $\omega$ and $-\omega$ give the same Frobenius manifold.
\end{remark}

\begin{definition}\label{def:preprim}
Let $(M,E,\nabla,R_\infty,\Phi,R_0,g)$ be a pre-Saito structure of weight~$w$. Let $\omega\in E$ be a $\nabla$-horizontal section. We will say that $\omega$ is \emph{pre-primitive} if it satisfies the following properties:
\begin{enumerate}
\item\label{def:preprim1}
$\omega^o$ and its images under the iterates of $\Phi^o_\xi$ ($\xi\in \Theta_M^o$) generate $E^o$,
\item\label{def:preprim2}
$\varphi_\omega^o:\Theta_M^o\to E^o$ is injective.
\end{enumerate}
We say that $\omega$ is \emph{strongly pre-primitive} if it satisfies moreover
\begin{enumerate}\setcounter{enumi}{2}
\item\label{def:preprim3}
$\omega\not\in\im\varphi_\omega$.
\end{enumerate}
\end{definition}

The third condition will only be useful when considering tensor products. Let us notice that, because of this condition, a primitive section is \emph{not} strongly pre-primitive. Let us also notice that the generating condition is somewhat stronger than what is needed to apply the results of Hertling and Manin, as $R_0^o$ is not used in \ref{def:preprim}\eqref{def:preprim1}. This will also be useful when considering tensor products. On the other hand, adding a new parameter as in Example \ref{exam:deformation} enables us to skip $R_0^o$ in the generating condition of Hertling and Manin.

From Corollary \ref{cor:hermanuniv} we get:

\begin{corollary}\label{cor:Frobuniv}
Let $(M,E,\nabla,R_\infty,\Phi,R_0,g)$ be a pre-Saito structure of weight~$w$. Let $\omega$ be a $\nabla$-horizontal pre-primitive section of $E$. Let us moreover assume that~$\omega$ is homogeneous with respect to $R_\infty$. Then, on the base space $N$ of any universal deformation of $(M,E,\nabla,R_\infty,\Phi,R_0,g)$ exists a canonical Frobenius manifold structure. The Frobenius manifold structures on two such deformations $N$ and $N'$ are isomorphic by an isomorphism which induces the identity on $M$.
\end{corollary}

\begin{proof}
The $\nabla$-horizontal extension $\omega_N$ of $\omega$ on any universal deformation space $N$ (whose existence is granted by Corollary \ref{cor:hermanuniv}) is primitive and remains homogeneous. It defines therefore a Frobenius manifold structure on $N$. Given another deformation with base space $N'$, it is obtained by pull-back by $\nu:N'\to N$. We have $E_{N'}=\cO_{N'}\otimes_{\cO_N}E_{N}$ and $\omega_{N'}=1\otimes\omega_N$. Keeping notation of the proof of Corollary \ref{cor:hermanuniv}, we have $\chi_{\omega_{N'}}=\chi_{\omega_N}\circ\nu$, hence $\varphi_{\omega_{N'}}=\varphi_{\omega_N}\circ T\nu$, and the structures on $\Theta_{N'}$ and $\Theta_N$ correspond through the isomorphism $T\nu$ (for the metric, one uses Corollary \ref{cor:hermanunivwith}).
\end{proof}

\begin{remark}
For $(M,E,\nabla,R_\infty,\Phi,R_0,g)$ and $\omega$ pre-primitive and homogeneous, there is a meaning to speak of \emph{the} Frobenius manifold structure determined by the pre-primitive homogeneous section $\omega$ on \emph{the} universal deformation.
\end{remark}

\begin{example}\label{exam:deformation2}
Let $(E^o,R_\infty,R_0^o,\omega^o)$ be as in Example \ref{exam:deformationgenerating} with $\rk E^o\geq2$. Then $\omega$ is strongly pre-primitive. Indeed, we have $\varphi_{\omega^o}(\partial_x)=R_0^o(\omega^o)\not\in\CC\cdot\omega^o$ as $\rk E^o\geq2$.

Assume moreover that $\omega^o$ and $R_0^o(\omega^o)$ are eigenvectors of $R_\infty$ with respective eigenvalues $\delta_0$ and $\delta_1=\delta_0-1$. Then the germ $\wt E^o=(E^o,\CC\cdot R_0^o(\omega^o))$ gets equipped with the structure of a Frobenius manifold. The Euler vector field $\gE$ is tangent to the line $M=\CC\cdot R_0^o(\omega^o)$ and, in the coordinates $x_0,\dots,x_{d-1}$ considered in \loccit, $\gE_{|M}=\partial_{x_1|M}$. The subsheaf of algebras $\cO_M[\gE_{|M}]\subset (\Theta_{\wt E^o|M},\star)$ is isomorphic to $\cO_M[y]/p(e^{-x_1}y)$, if $p$ denotes the characteristic polynomial of $R_0^o$ and the inclusion above is in fact an equality.
\end{example}

\begin{example}[Quantum cohomology of the projective space]\label{exam:proj}
Let us consider the pre-Saito structure $(E^o,R_\infty,R_0^o,g^o)$ equipped with the pre-primitive form $\omega^o$ given by the following data:
\begin{itemize}
\item
$E^o$ is $\CC^{n+1}$ with its canonical basis $\omega^o=\omega_0^o,\omega^o_1,\dots,\omega^o_n$,
\item
the matrix of $R_0^o$ is
\[
(n+1)\begin{pmatrix}
0&\cdots&\cdots&0&1\\
1&0&\cdots&\cdots&0\\
0&1&0&\cdots&0\\
\vdots&\ddots&\ddots&\ddots&\vdots\\
0&\cdots&0&1&0
\end{pmatrix}
\]
and that of $R_\infty$ is $-\diag(0,1,\dots,n)$,
\item
we have $g^o(\omega^o_k,\omega^o_\ell)=1$ if $k+\ell=n$ and $0$ otherwise.
\end{itemize}
The germ of universal Frobenius manifold defined by $(E^o,R_\infty,R_0^o,g^o,\omega^o)$ is equal to that defined by the quantum cohomology of $\PP^n$ (\cf \cite[\S II.4]{Manin96}). Let us denote by $t_0,\dots,t_n$ the flat coordinates corresponding to the basis $\omega_0^o,\dots,\omega_n^o$.

The trivial deformation parametrized by $\Afu$ is given by the linear map $\chi_\omega(x)=(n+1)x\omega^o_1$, and the Frobenius manifold structure is defined along this line. Working in the flat coordinate $t_1=(n+1)x$, the pre-Saito structure at the point $t_1$ is $(E^o,R_\infty,R_0(t_1),g^o)$, with
\[
R_0(t_1)=(n+1)\begin{pmatrix}
0&\cdots&\cdots&0&e^{t_1}\\
1&0&\cdots&\cdots&0\\
0&1&0&\cdots&0\\
\vdots&\ddots&\ddots&\ddots&\vdots\\
0&\cdots&0&1&0
\end{pmatrix}=-(n+1)\Phi_{\partial_{t_1}}(t_1).
\]
Moreover, by construction, we have $\Phi_{\partial_{t_i}}(\omega^o)_{|\Afu}=-\omega_i^o$ for any $i$, and therefore, denoting now $\Phi_{\partial_{t_i}}=-\partial_{t_i}\star{}$, we get
\[
(\partial_{t_i}\star{})_{|\Afu}=(\partial_{t_1}\star{})^i_{|\Afu}\qquad i=0,\dots,n.
\]
In this example, $\wt E^o$ denotes the germ of $E^o$ along $\Afu=\CC\cdot\omega_1^o$, equipped with the flat coordinates $(t_0,\dots,t_n)$, and we have an isomorphism of sheaves of algebras $(\Theta_{\wt E^o|\Afu},\star_{|\Afu})\isom\cO_{\Afu}[y]/(p(t_1,y))$, with $p(y)=y^{n+1}-e^{t_1}$, given by $\partial_{t_k}\mto[y^k]$.

Let us also note that the Frobenius manifold structure on $\wt E^o$ is invariant by translation of $t_1$ by $2i\pi\ZZ$.
\end{example}

\section{Application to tensor products and alternate products}\label{sec:tensor}

\subsection{Tensor product of two Frobenius manifolds}\label{subsec:tensor}
Let us show how the previous results enable us to recover existence and uniqueness results concerning tensor products of Frobenius manifolds (see \cite{Kaufmann99} and also \cite[\S III.7]{Manin96}).

Let us start with pre-Saito structures. The tensor product of two pre-Saito structures $(M,E,\nabla,R_\infty,\Phi,R_0,g)$ and $(M,E',\nabla',R'_\infty,\Phi',R'_0,g')$ on a given manifold $M$, of respective weights $w$ and $w'$ is $(M,E'',\nabla'',R''_\infty,\Phi'',R''_0,g'')$, with
\begin{itemize}
\item
$E''=E\otimes_{\cO_M}E'$,
\item
$\nabla''=\nabla\otimes\id+\id\otimes\nabla'$,
\item
$R''_\infty=R_\infty\otimes\id+\id\otimes R'_\infty$,
\item
$\Phi''=\Phi\otimes\id+\id\otimes\Phi'$,
\item
$R''_0=R_0\otimes\id+\id\otimes R'_0$,
\item
$g''(e\otimes e',f\otimes f')=g(e,f)g'(e',f')$.
\end{itemize}
This produces a pre-Saito structure of weight $w''=w+w'$. Let us note that, from the point of view of the connection $\nablab$ defined in \eqref{eq:nablab}, the tensor product is associated to the tensor product connection on $\bE\otimes_{\cO_M[\hb]}\bE'$. Given $r\geq2$, we can similarly define the $r$-fold tensor product, the $r$-fold symmetric product and the $r$-fold alternate product of a pre-Saito structure. From the point of view of $(\bE,\nablab)$, they correspond respectively to the natural connection $\nablab$ on
\[
\otimes^r_{\cO_M[\hb]}\bE,\quad \Sym^r_{\cO_M[\hb]}\bE,\quad \wedge^r_{\cO_M[\hb]}\bE.
\]

Let us now consider two pre-Saito structures $(M_i,E_i,\nabla_i,R_{\infty,i},\Phi_i,R_{0,i},g_i)$ of weights~$w_i$ (\hbox{$i=1,2$}). Let us denote by $p_1,p_2$ the projections $M_1\times M_2\to M_1,M_2$. The \emph{external tensor product} of these pre-Saito structures is, by definition,
\[
p_1^*(M_1,E_1,\nabla_1,R_{\infty,1},\Phi_1,R_{0,1},g_1)\otimes_{\cO_{M_1\times M_2}}p_2^*(M_2,E_2,\nabla_2,R_{\infty,2},\Phi_2,R_{0,2},g_2),
\]
where the pull-back $p_i^*$ ($i=1,2$) has been defined in \S\ref{subsec:univdef}. The external tensor product has weight $w_1+w_2$. We will denote it by $\boxtimes$.

\begin{lemma}\label{lem:tensorpreprim}
Assume that $\omega_1,\omega_2$ are strongly pre-primitive horizontal sections of $E_1,E_2$. Then $\omega=\omega_1\boxtimes\omega_2\in E_1\boxtimes E_2$ is a strongly pre-primitive horizontal section of $E_1\boxtimes E_2$.

If moreover $\omega_1,\omega_2$ are homogeneous of respective degrees $q_1,q_2$, then $\omega$ is homogeneous of degree $q_1+q_2$.
\end{lemma}

\begin{proof}
If we denote by $x$ the coordinates on $M_1$ and by $y$ that on $M_2$, we remark that $\Phi_{\partial_{x_i}}(\omega)=(\Phi_{1,\partial_{x_i}}(\omega_1))\boxtimes\omega_2$ and a similar result for $\Phi_{\partial_{y_j}}(\omega)$. Therefore, $\omega$ and the iterates of $\Phi_\xi$ ($\xi\in\Theta_{M_1\times M_2}$) acting on $\omega$ generate $E_1\boxtimes E_2$.

Moreover, $\omega_1^o\otimes\omega_2^o$ does not belong to the vector space generated by the $\varphi_{\omega_1^o}(\partial_{x_i})\otimes\nobreak\omega_2^o$ and the $\omega_1^o\otimes\varphi_{\omega_2^o}(\partial_{y_j})$, which clearly form part of a basis of $E_1^o\otimes E_2^o$, hence the strong pre-primitivity.

Lastly, the homogeneity condition for $\omega$ directly follows from the formulas above.
\end{proof}

The reason to impose the third condition in the definition of a strongly pre-primitive section is to insure that, in the previous lemma, $\varphi_\omega$ remains injective. Otherwise, if we set $e_i=\varphi_{\omega_i}^{-1}(\omega_i)\in\Theta_{M_i}$ ($i=\nobreak1,2$), and if we denote similarly the corresponding vector field on $M_1\times M_2$, then $\varphi_\omega(e_1-\nobreak e_2)=\omega_1\boxtimes\omega_2-\omega_1\boxtimes\omega_2=0$, so $\varphi_\omega$ is not injective. Let us also notice that the lemma holds if only one of both pre-primitive sections $\omega_1$ and $\omega_2$ is strong.

In conclusion, \emph{the tensor product is well-defined for pre-Saito structures equipped with a strongly pre-primitive homogeneous section}. We will say that the Frobenius manifold structure associated according to Corollary \ref{cor:Frobuniv} to this tensor product is the \emph{tensor product} of the Frobenius manifold structures corresponding to each term, although this is incorrect, strictly speaking. (Another approach of the tensor product is given in \cite{Kaufmann99}, see also \cite[\S III.7]{Manin96}).

\begin{example}\label{exam:deformation3}
Let $(E^o,R_\infty,R^o_0)$ and $(E^{\prime o},R'_\infty,R^{\prime o}_0)$ be two pre-Saito structures (without metric), with underlying manifold $M,M'$ reduced to a point. We define their tensor product $(E^{\prime\prime o},R''_\infty,R^{\prime\prime o}_0)$ as an object of the same kind:
\begin{align*}
E^{\prime\prime o}&=E^o\otimes_\CC E^{\prime o},\\
R''_\infty& =R_\infty\otimes\id+\id\otimes R'_\infty,\\
R^{\prime\prime o}_0&=R^o_0\otimes \id +\id\otimes R^{\prime o}_0.
\end{align*}

Assume that there exist $\omega^o\in E^o$ and $\omega^{\prime o}\in E^{\prime o}$ such that the $(R^o_0)^k(\omega^o)$ ($k\geq0$) generate $E^o$, and similarly with ``prime''. Then, by Corollary \ref{cor:hermanuniv}, there exists a universal deformation of $(E^o,R_\infty,R^o_0)$ and $(E^{\prime o},R'_\infty,R^{\prime o}_0)$. However, the \hbox{$(R^o_0\otimes \id+\id\otimes R^{\prime o}_0)^k(\omega^o\otimes\omega^{\prime o})$} may not generate $E^o\otimes_\CC E^{\prime o}$, and the same corollary cannot be applied to the tensor product.

We can use Example \ref{exam:deformation} to overcome this difficulty. Indeed, let us denote by $(M,E,\nabla,R_\infty,\Phi,R_0)$ and $(M',E',\nabla',R'_\infty,\Phi',R'_0)$ the one-parameter deformations of $(E^o,R_\infty,R^o_0)$ and $(E^{\prime o},R'_\infty,R^{\prime o}_0)$ defined there. The external tensor product is defined as above on $M''=M\times M'$ and $E''=p^*E\otimes_{\cO_{M''}}p^{\prime*}E'$ ($p,p'$ the projections from $M''$ to $M,M'$), adding the relations
\[
\nabla''=\nabla\otimes\id+\id\otimes\nabla',\quad
\Phi''=\Phi\otimes\id+\id\otimes\Phi''.
\]
We thus have
\[
\Phi''_{\partial_x|x=x'=0}=-R_0^o\otimes\id,\quad
\Phi''_{\partial_{x'|x=x'=0}}=-\id\otimes R_0^{\prime o}.
\]
Let us note that the flat extensions of $\omega^o$ and $\omega^{\prime o}$ are now strongly pre-primitive. From Lemma \ref{lem:tensorpreprim}, we conclude that the flat extension $\omega''$ of $\omega^o\otimes\omega^{\prime o}$ is strongly pre-primitive, and therefore the generating condition of Proposition \ref{prop:herman} is fulfilled (even without using $R''_0$), so a universal deformation of $(M'',E'',\nabla'',\Phi'',R_0'')$ does exist. Moreover, according to Example \ref{exam:deformation2} and if metrics $g^o$, $g^{\prime o}$ do exist, giving weights $w,w'$, the tensor product of the corresponding Frobenius manifold structures is well-defined if we moreover assume that $\omega^o,\omega^{\prime o}$ are homogeneous.
\end{example}

\subsection{Symmetric and alternate product of a Frobenius manifold}\label{subsec:alternate}
Let us fix a pre-Saito structure $(M,E,\nabla,R_\infty,\Phi,R_0,g)$ of weight~$w$. For any $r\geq1$, we can consider the $r$-fold external tensor product as in \S\ref{subsec:tensor}, with base space $M^r$ and vector bundle ${\boxtimes}^rE$. Assume that $\omega$ is a strongly pre-primitive (\resp homogeneous) flat section of $E$. Then, we have seen that ${\boxtimes}^r\omega$ is so for ${\boxtimes}^rE$.

Moreover, we have a natural action of the symmetric group $\rW=\gS_r$ on the pre-Saito structure ${\boxtimes}^r(M,E,\nabla,R_\infty,\Phi,R_0,g)$, and ${\boxtimes}^r\omega$ is $\rW$-invariant.

It follows from Corollary \ref{cor:Frobuniv} that $\wt{\otimes^rE^o}$ is equipped with the $r$-fold tensor product Frobenius manifold structure, and the natural action of $\rW$ is compatible with this structure.

From now on, we will only consider the situation of Example \ref{exam:deformation2} (in particular, we assume that $\omega^o$ and $R^o_0(\omega^o)$ below are eigenvectors of $R_\infty$ with respective eigenvalues $\delta_0$ and $\delta_1=\delta_0-1$). Let us fix a pre-Saito structure $(E^o,R_\infty,R^o_0,g^o)$ of weight $w$ (with base manifold $M$ reduced to a point) with a homogeneous $R_0^o$-cyclic vector $\omega^o\in E^o$. Let $r$ be an integer $\geq2$. Now, there is no difference between the external tensor product and the tensor product over~$\CC$. On $\otimes^rE^o$ we have operators denoted by $\oR_\infty,\oR_0^o$: for instance,
\[
\oR_0^o =\sum_{i=1}^r\id\otimes\cdots\otimes\id\otimes\underset{i}{R_0^o}\otimes\id\otimes\cdots\otimes\id.
\]
These operators are $\rW$-invariant. Therefore, they induce on the symmetric product $\Sym^rE^o\defin(\otimes^rE^o)^\rW$ and on the alternate product $\wedge^rE^o$ similar operators.

Let $(\Afu,E,\nabla,R_\infty,\Phi,R_0)$ be the one-parameter deformation of $(E^o,R_\infty,R^o_0)$ constructed in Example \ref{exam:deformation}, together with the flat extension $\omega=1\otimes\omega^o$ of~$\omega^o$. Then, by assumption on $\omega^o$ and according to Example \ref{exam:deformation2}, $\omega$ is strongly pre-primitive and homogeneous. 

\subsubsection*{The $r$-fold tensor product}
The external tensor product ${\boxtimes}^r(\Afu,E,\nabla,R_\infty,\Phi,R_0)$ is a $r$-parameter deformation of $(\otimes^rE^o,\oR_\infty,\oR_0^o)$, equipped with the strongly pre-primitive homogeneous section ${\boxtimes}^r\omega=1\otimes(\otimes^r\omega^o)$. Its germ at the origin has a universal deformation with base manifold equal to the germ $\wt{\otimes^rE^o}$ of $\otimes^rE^o$ at $0$ that we denote by
\begin{equation}\label{eq:univprefro}
(\wt{\otimes^rE^o},\cO_{\wt{\otimes^rE^o}}\otimes_\CC (\otimes^rE^o),d,\oR_\infty,\wt\Phi,\wt R_0).
\end{equation}
The tangent map $\varphi_{{\scriptscriptstyle\boxtimes}^r\omega}$ of the embedding $\chi_{{\scriptscriptstyle\boxtimes}^r\omega}:((\Afu)^r,0)\hto\wt{\otimes^rE^o}$ sends the $j$-th vector basis of $\Theta^o_{(\Afu)^r}$ to
\[
\omega^o\otimes\cdots\otimes\omega^o\otimes \underset{j}{R_0^o(\omega^o)}\otimes\omega^o\otimes\cdots\otimes\omega^o.
\]
This universal deformation defines a Frobenius manifold structure $(\wt{\otimes^rE^o},\star,\otimes^rg,e,\gE)$ of weight $D=2-rw$. The natural action of $\rW$ on $\wt{\otimes^rE^o}$ is by automorphisms of the Frobenius manifold structure.

Let us set $d=\dim_\CC E^o$. For any multi-index $\alpha\in\{0,\dots,d-1\}^r$ we set $e_\alpha^o=(R_0^o)^{\alpha_1}(\omega^o)\otimes\cdots\otimes(R_0^o)^{\alpha_r}(\omega^o)$, getting thus a basis $\bme^o$ of $\otimes^rE^o$. We denote by $(x_\alpha)$ the corresponding coordinates on $\wt{\otimes^rE^o}$.

\begin{lemma}\label{lem:Phialpha}
For any multi-index $\alpha$, we have
\[
\wt\Phi_{\partial_{x_\alpha}}^o=-(R_0^o)^{\alpha_1} \otimes\cdots\otimes(R_0^o)^{\alpha_r}.
\]
\end{lemma}

\begin{proof}
If $1\otimes(\otimes^r\omega^o)$ denotes the horizontal extension of $\otimes^r\omega^o$ on $\wt{\otimes^rE^o}$, we have by definition $\varphi_{1\otimes(\otimes^r\omega^o)}^o(\partial_{x_\alpha})=e_\alpha^o$, that is, $\wt\Phi_{\partial_{x_\alpha}}^o(\otimes^r\omega^o)=-e^o_\alpha$. On the other hand, the images of $\otimes^r\omega^o$ under the iteration of the operators $\wt\Phi_{\partial_{x_j}}^o$, with $j=1,\dots,r$, generate $\otimes^rE^o$. The commutation relations $[\wt\Phi_{\partial_{x_\alpha}}^o,\wt\Phi_{\partial_{x_j}}^o]=0$ imply that $\wt\Phi_{\partial_{x_\alpha}}^o$ is determined by its value on $\otimes^r\omega^o$, hence the assertion.
\end{proof}

\begin{remark}\label{rem:A1r}
The previous results hold all along $(\Afu)^r\subset \otimes^rE^o$ and not only at the origin, so that $\wt{\otimes^rE^o}$ can be regarded as the analytic germ of $\otimes^rE^o$ along $(\Afu)^r$.
\end{remark}

\subsubsection*{The $r$-fold symmetric product}
The space $(\otimes^rE^o)^\rW=\Sym^rE^o$ has a basis obtained by symmetrization of the basis $\bme^o$ of $\otimes^rE^o$. We will consider the subspace $\Elem^rE^o$ generated by the symmetrization of the vectors $e^o_\alpha$ with $\alpha_j\in\{0,1\}$ for any $j$ and $\alpha_k=1$ for at least one $k$. It has dimension $r$.

\begin{lemma}\label{lem:elemsym}
The restriction
\[
\big(\wt{\Elem^rE^o},\cO_{\wt{\Elem^rE^o}}\otimes_\CC (\otimes^rE^o)^\rW,d,\oR_\infty,\wt\Phi,\wt R_0\big)
\]
is a pre-Saito structure admitting $1\otimes(\otimes^r\omega^o)$ as pre-primitive homogeneous section.
\end{lemma}

\begin{proof}
A basis of the tangent space to $\wt{\Elem^rE^o}$ at $o$, that is, $\Elem^rE^o$, consists of the elementary symmetric vector fields $\xi_k=\sum_{\alpha\in A_k}\partial_{x_\alpha}$, where $A_k=\{\alpha\in\{0,1\}^r\mid\sum_j\alpha_j=k\}$, and $k=1,\dots,r$. Any element of $(\otimes^rE^o)^\rW$ can be obtained from $\otimes^r\omega^o$ by applying a symmetric polynomial in the $\wt\Phi_{\partial_{x_\alpha}}^o$. By the lemma above, it can thus be obtained by applying iterations of the elementary operators $\wt\Phi_{\xi_k}^o$.
\end{proof}

Lemma \ref{lem:elemsym}, together with Corollary \ref{cor:Frobuniv}, endows $\wt{\Sym^rE^o}$ with the structure of a Frobenius manifold, through the infinitesimal period mapping defined by $1\otimes(\otimes^r\omega^o)$.

On the other hand, let us consider the restriction of \eqref{eq:univprefro} to $\wt{\Sym^rE^o}$ or to $\wt{\Elem^rE^o}$. The action of $\rW$ on the base manifold is equal to the identity, so these restrictions have a linearized $\rW$-action. We claim that
\[
\big(\wt{\Sym^rE^o},\cO_{\wt{\Sym^rE^o}}\otimes_\CC (\otimes^rE^o),d,\oR_\infty,\wt\Phi,\wt R_0\big)
\]
is the universal deformation with linearized $\rW$-action of
\[
\big(\wt{\Elem^rE^o},\cO_{\wt{\Elem^rE^o}}\otimes_\CC (\otimes^rE^o),d,\oR_\infty,\wt\Phi,\wt R_0\big).
\]
This follows from the $\rW$-invariant version of Proposition \ref{prop:herman} and Corollary \ref{cor:hermanuniv} explained after Corollary \ref{cor:hermanWuniv}.

\begin{remark}\label{rem:diag}
As in Remark \ref{rem:A1r}, the results above hold all along the $\rW$-invariant part of $(\Afu)^r$, that is, the diagonal $\Afu$, and we can regard $\wt{\Elem^rE^o}$ or $\wt{\Sym^rE^o}$ as the analytic germs of $\Elem^rE^o$ or $\Sym^rE^o$ along the diagonal $\Afu$.
\end{remark}

\subsubsection*{The $r$-fold alternate product}
We assume here that $r<d$. In order to obtain a Frobenius manifold structure on the subvariety $\wt{\wedge^rE^o}$, we do not use the same procedure as for $\wt{\Sym^rE^o}$, as $\rW$ does not act trivially on this subvariety. We notice however that the bundle $\cO_{\wt{\Sym^rE^o}}\otimes_\CC (\otimes^rE^o)$, and therefore $\cO_{\wt{\Elem^rE^o}}\otimes_\CC (\otimes^rE^o)$, is equipped with a linearized $\rW$-action. We can thus consider the anti-invariant subbundle $\cO_{\wt{\Sym^rE^o}}\otimes_\CC (\wedge^rE^o)$, which is left invariant by $\oR_\infty$, $\wt R_0$ and $\wt\Phi_\xi$ for any vector field $\xi$ tangent to $\wt{\Sym^rE^o}$.

On $\wt{\Sym^rE^o}$ (hence on $\wt{\Elem^rE^o}$) exists the anti-invariant part of the restriction of the pre-Saito structure \eqref{eq:univprefro} to $\wt{\Sym^rE^o}$, which is a pre-Saito structure that we denote by
\begin{equation}\label{eq:antiinv}
\big(\wt{\Sym^rE^o},\cO_{\wt{\Sym^rE^o}}\otimes_\CC (\wedge^rE^o),d,\oR_\infty,\wt\Phi,\wt R_0\big).
\end{equation}

\begin{lemma}
The restriction
\[
\big(\wt{\Elem^rE^o},\cO_{\wt{\Elem^rE^o}}\otimes_\CC (\wedge^rE^o),d,\oR_\infty,\wt\Phi,\wt R_0\big)
\]
is a pre-Saito structure admitting
\[
1\otimes\wt\omega^o\defin1\otimes\big(\omega^o\wedge R_0^o(\omega^o)\wedge\cdots\wedge(R_0^o)^{r-1}(\omega^o)\big)
\]
as pre-primitive section. If $R_0^o(\omega^o),\dots,(R_0^o)^{r-1}(\omega^o)$ are eigenvectors of $R_\infty$, then $\wt\omega^o$ is homogeneous.
\end{lemma}

\begin{proof}
The homogeneity condition is clear, according to the assumption. Let us also note that, for $k=1,\dots,r$,
\[
\wt\Phi_{\xi_k}^o[\wt\omega^o]=-\omega^o\wedge\cdots\wedge (R_0^o)^{r-k-1}(\omega^o)\wedge(R_0^o)^{r-k+1}(\omega^o)\wedge\cdots\wedge(R_0^o)^{r}(\omega^o),
\]
so the injectivity condition is clear.

In order to check the generating condition, it is convenient to use the presentation of the algebra $(E^o,\star)$ as $\CC[y]/p(y)$, where $y$ denotes $R_0^o(\omega^o)$ and $p$ is the minimal polynomial of $R_0^o$. Then $\otimes^r(E^o,\star)=\CC[y_1,\dots,y_r]/(p(y_1),\dots,p(y_r))$. Any \hbox{(anti-)}invariant element of $\otimes^rE^o$ has a representative in $\CC[y_1,\dots,y_r]$ which is \hbox{(anti-)}invariant (by taking the \hbox{(anti-)}symmetrization of any representative), and $\wt\omega^o=[1]\wedge[y]\wedge\cdots\wedge[y^{r-1}]$ is the class of $\prod_{i>j}(y_i-y_j)$. Moreover, it is easy to check that $\prod_{i>j}(y_i-y_j):\CC[y_1,\dots,y_r]^\rW\to\CC[y_1,\dots,y_r]^\ant$ is onto. On the other hand, $\CC[y_1,\dots,y_r]^\rW$ is generated by the elementary symmetric polynomials. This gives the generating condition for $\wt\omega^o$.
\end{proof}

\begin{corollary}\label{cor:Frobuniv2}
Let $(E^o,R_\infty,R_0^o,g^o)$ be a punctual pre-Saito structure with $\dim E^o=d$. Let $\omega^o$ be a cyclic vector for $R_0^o$. Assume that $\omega^o,\dots,(R_0^o)^{r-1}(\omega^o)$ are eigenvectors of $R_\infty$ and, as above, that $\delta_1=\delta_0-1$. Then the $r$-fold alternate product of the (germ of) Frobenius manifold $\wt E^o$ is well-defined as the Frobenius manifold attached to the universal deformation of the pre-Saito structure
\[
\big(\wt{\Elem^rE^o},\cO_{\wt{\Elem^rE^o}}\otimes_\CC (\wedge^rE^o),d,\oR_\infty,\wt\Phi,\wt R_0,\otimes^rg^o\big),
\]
with primitive homogeneous section $\wt\omega^o$.\qed
\end{corollary}

Let us denote by $\star$ the product on $\Theta_{\wt{\otimes^rE^o}}$ given by the Frobenius manifold structure constructed above. By Lemma \ref{lem:Phialpha}, we have
\[
\varphi_{\otimes^r\omega^o}(\partial_{x_\alpha}\star\partial_{x_\beta})=(R_0^o)^{\alpha_1+\beta_1}(\omega^o) \otimes\cdots\otimes(R_0^o)^{\alpha_r+\beta_r}(\omega^o),
\]
where, if $\alpha_j+\beta_j\geq d$, we expand $(R_0^o)^{\alpha_j+\beta_j}(\omega^o)$ in terms of the $(R_0^o)^k(\omega^o)$, with $k=0,\dots,d-1$. Then, from Corollary \ref{cor:Frobuniv} we can give a realization of the $r$-fold alternate Frobenius structure:

\begin{corollary}\label{cor:Frobuniv3}
Under the assumptions of Corollary \ref{cor:Frobuniv2}, let $N$ be a germ of complex manifold with $\wt{\Elem^rE^o}\subset N\subset\wt{\Sym^rE^o}$ such that the product with $\wt\omega^o$:
\[
{}\star\wt\omega^o:\Theta^o_{\wt{\Sym^rE^o}}=\Sym^rE^o\to \wedge^rE^o
\]
induces an isomorphism $\Theta^o_N\to\wedge^rE^o$. Then the restriction to $N$ of \eqref{eq:antiinv} is a universal deformation of its restriction to $\wt{\Elem^rE^o}$, and the primitive homogeneous section $1\otimes\wt\omega^o$ induces a Frobenius manifold structure on $N$, which is independent, up to isomorphism, on the choice of $N$.\qed
\end{corollary}

\begin{remark}\label{rem:diag2}
As in Remark \ref{rem:diag}, one can notice that Corollary \ref{cor:Frobuniv2} holds all along the diagonal $\Afu$ and that Corollary \ref{cor:Frobuniv3} holds on any open set of the diagonal on which the isomorphism condition on ${}\star\wt\omega^o$ is satisfied.
\end{remark}

\subsection{Quantum cohomology of the Grassmannian as an alternate product of a Frobenius manifold}
In this paragraph, we consider Example \ref{exam:proj} with its notation and we take $r\leq n$. The assumptions of Corollary \ref{cor:Frobuniv2} are then satisfied. The germ $\wt{\Elem^rE^o}$ is now a germ along the diagonal $\Afu=\CC\cdot(\omega^o_1\otimes\cdots\otimes\omega^o_1)\subset\wt{\Sym^rE^o}$. We can apply the previous results all along the diagonal. We also replace above $R_0^o$ with $\partial_{t_1}\star$ and we set $\wt\omega^o=\omega_0^o\wedge\omega_1^o\wedge\cdots\wedge\omega_{r-1}^o$. For any multi-index $\alpha\in\{0,\dots,n\}^r$, we also set $\omega_\alpha^o=\omega_{\alpha_1}^o\otimes\cdots\otimes\omega_{\alpha_r}^o$ and we denote by $\partial_{t_\alpha}$ the corresponding germs of vector fields on $\wt{\otimes^rE^o}$ along $(\Afu)^r$. We will denote by $\bun_i$ the multi-index $\alpha$ with $\alpha_j=\delta_{ij}$ for all $j=1,\dots,r$. The infinitesimal period mapping induces an isomorphism of algebras
\begin{equation}
\tag{$*$}
\begin{split}
(\Theta_{\wt{\otimes^rE^o}},\star)_{|(\Afu)^r}&\to\cO_{(\Afu)^r}[y_1,\dots,y_r]/(p(t_{\bun_1},y_1),\dots,p(t_{\bun_r},y_r)),\\
(\partial_{t_\alpha})_{|(\Afu)^r}&\mto [y^\alpha],
\end{split}
\end{equation}
with $p(t,y)=y^{n+1}-e^{t}$.

Along the diagonal $\Afu$ with coordinate $t$, the tangent algebra of the $r$-fold alternate product Frobenius manifold is $\big[\cO_{\Afu}[y_1,\dots,y_r]/(p(t,y_1),\dots,p(t,y_r))\big]^\ant$ and $1\otimes\wt\omega^o$ is the class of $\prod_{i>j}(y_i-\nobreak y_j)$. Let us also remark that this Frobenius structure is invariant under the translation of $t$ by $2i\pi\ZZ$.

From the main result in \cite{CF-K-S06} we obtain:
\begin{theorem}\label{th:cfks}
The Frobenius manifold structure attached to the quantum cohomology of the complex Grassmannian $G(r,n+1)$ of $r$-planes in $\CC^{n+1}$ is isomorphic to the germ at $t^o=(r-1)i\pi\in\Afu$ of the $r$-fold alternate Frobenius manifold structure of the quantum cohomology of $\PP^n$ defined through the pre-primitive homogeneous section $\wt\omega\defin\rho_r(1\otimes\wt\omega^o)$, with $\rho_r^2=(-1)^{\binom r2}/r!$.
\end{theorem}

Let us note that the choice of a square root of $(-1)^{\binom r2}/r!$ is not important, according to Remark \ref{rem:sqrt}.

\begin{proof}
Let us denote by $\PP$ the $r$-fold product of $\PP^n$. We set $E^o=H^*(\PP^n)$ and $(\omega_i^o)_{i=0,\dots,n}$ is the basis generated by the hyperplane class $H=\omega_1^o$. Then $\otimes^rE^o=H^*(\PP)$. The Frobenius structure attached to the quantum cohomology of $\PP$ is known to be the $r$-fold tensor product of that of $\PP^n$.

If $S(Y_1,\dots,Y_r)$ is any polynomial in $r$ variables with degrees in each variable belonging to $[0,n]$, we define $S(\omega_0^o\otimes\cdots\otimes\omega_0^o)\in \otimes^rE^o$ by replacing each monomial $Y_1^{m_1}\cdots Y_r^{m_r}$ in $S$ by $\omega_{m_1}^o\otimes\nobreak\cdots\otimes\nobreak\omega_{m_r}^o$.
Let us denote by $s_\lambda$ the Schur polynomials in $r$ variables indexed by partitions $\lambda$ having Young diagrams in a rectangle $r\times(n+1-r)$ and let $N\subset \Sym^rE^o$ be the linear subspace having $\omega^o_{s_\lambda}\defin\nobreak s_\lambda(\omega_0^o\otimes\nobreak\cdots\otimes\nobreak\omega_0^o)$ as a basis. In particular, $N\supset\Elem^rE^o$ and the ordinary cup product $\wt\omega^o\cup{}:\Sym^rE^o=H^*(\PP)^\rW\to H^*(\PP)^\ant=\wedge^rE^o$ induces an isomorphism $N\isom H^*(\PP)^\ant$.

Let us denote by $\xi_{s_\lambda}$ the vector field on $N$ corresponding to $\omega^o_{s_\lambda}$. By definition of the product $\star$ and the isomorphism $(*)$ above, and as $\xi_{s_\lambda}$ is a linear combination of the $\partial_{t_\alpha}$, the restriction $\xi_{s_\lambda|(\Afu)^r}$ is sent to $[s_\lambda(y_1,\dots,y_r)]$. Proving that the isomorphism condition of Corollary \ref{cor:Frobuniv3} is satisfied for~$N$ amounts then to proving that
\begin{multline*}
\prod_{i>j}(y_i-y_j):[\cO_{\Afu}[y_1,\dots,y_r]/(p(t,y_1),\dots,p(t,y_r))]^\rW\\
\to [\cO_{\Afu}[y_1,\dots,y_r]/(p(t,y_1),\dots,p(t,y_r))]^\ant
\end{multline*}
induces an isomorphism on the subsheaf generated by the $[s_\lambda(y_1,\dots,y_r)]$.

The sheaf $\cO_{\Afu}[y_1,\dots,y_r]/(p(t,y_1),\dots,p(t,y_r))$ is filtered according to the total degree in $y_1,\dots,y_r$ and the graded sheaf is $\cO_{\Afu}[y_1,\dots,y_r]/(p_0(y_1),\dots,p_0(y_r))$ with $p_0(y)=y^{n+1}$. As the action of $\rW$ on $\cO_{\Afu}[y_1,\dots,y_r]/(p(t,y_1),\dots,p(t,y_r))$ strictly preserves the filtration, taking the (anti-)invariant subsheaf commutes with gradation. Lastly, the morphism induced by the multiplication by $\prod_{i>j}(y_i-y_j)$ induces the same morphism at the graded level.

Now, at the graded level, we recover the ordinary cup product $\wt\omega^o\cup{}:N\to \wedge^rE^o$, which is an isomorphism. We conclude that $\wt\omega^o\star{}:\Theta_{N|\Afu}\to\cO_{\Afu}\otimes_\CC(\wedge^rE^o)$ is an isomorphism, and Corollary \ref{cor:Frobuniv3} equips the germ of $N$ along $\Afu$ of a canonical Frobenius structure isomorphic to the $r$-fold alternate product of that attached to the quantum cohomology of $\PP^n$.

Now, Corollary \ref{cor:Frobuniv3}, when applied to the germ of $N$ at $t^o=(r-1)i\pi\in\Afu$, gives a Frobenius structure isomorphic to that attached to the cohomology of the Grassmannian $G(r,n+1)$ (up to the normalizing factor $\rho_r$): indeed, the main result in \cite{CF-K-S06} gives a similar statement, but working with a Novikov variable $Q$ (the sign change in the Novikov variables in \loccit\ amounts here to the translation of the variable $t$ by $(r-1)i\pi$); from the previous considerations, we conclude in particular that the Gromov-Witten potential on $N$ is convergent, hence we can set the Novikov variable to $1$ in the result of \loccit
\end{proof}

\section{Alternate Thom-Sebastiani}\label{sec:TS}

This section, which is independent of the previous one, gives the necessary tools for the geometric interpretation given in \S\ref{sec:geom} of the alternate product of Frobenius manifolds constructed in \S\ref{subsec:alternate}.

\subsection{Holonomic $\cD$-modules and perverse sheaves with an action of a~finite group}

Let $Z$ be a complex manifold (\resp a smooth complex algebraic variety) and let $\cD_Z$ be the sheaf of holomorphic (\resp algebraic) differential operators on~$Z$.

Let $\rW$ be a finite group equipped with a non trivial character $\sgn:\rW\to\{\pm1\}$. For instance, $\rW=\gS_r$ is the symmetric group on $r$ letters and $\sgn$ is the signature.

Let $\cM$ be a holonomic (left or right) $\cD_Z$-module equipped with an action of the group $\rW$ by $\cD_Z$-automorphisms. Let $\cM^\ant$ be the biggest submodule on which any $\rw$ in $\rW$ acts by $\sgn(\rw)$. Let $a_{\cM}:\cM\to\cM$ be the antisymmetrization map
\[
m\mto \frac{1}{|\rW|}\sum_{\rw\in\rW}\sgn(\rw)\rw(m).
\]
We denote by $D\cM$ the dual holonomic $\cD_Z$-module. It comes naturally equipped with a dual action of $\rW$.

\begin{proposition}\label{prop:antisym}
We have $\cM^\ant=\im a_{\cM}$ and a decomposition $\cM=\ker a_{\cM}\oplus\cM^\ant$. Moreover, we have $Da_{\cM}=a_{D\cM}$ and an isomorphism $D(\cM^\ant)\simeq(D\cM)^\ant$.
\end{proposition}

\begin{proof}
The first point follows from the identity $a_{\cM}\circ a_{\cM}=a_{\cM}$ and the identification $\cM^\ant=\ker(a_{\cM}-\id)$. That $Da_\cM=a_{D\cM}$ follows from the exactness of the contravariant functor $D$ on holonomic modules. The second assertion is then clear.
\end{proof}

\begin{remark}[$\QQ$-perverse sheaves]\label{rem:antisym}
The same result holds for $\QQ$-perverse sheaves~$\cG$ on any reduced analytic space $Z$, if $\cG$ is equipped with an action of $\rW$ by automorphisms, where $\ker,\coker$ and ${}^\ant$ are taken in the abelian category of $\QQ$-perverse sheaves. The point is that the antisymmetrization morphism $a_\cG$ is well-defined as a morphism in this category, as $\Hom_{\Perv(Z)}(\cG,\cG)$ is a $\QQ$-vector space.
\end{remark}

\begin{proposition}\label{prop:imdirant}
Let $g:Z\to Z'$ be a proper map (between complex analytic manifolds or between smooth complex algebraic varieties). If $\cM$ is as above, then for any $k\in\ZZ$, the $\cD_{Z'}$-modules $\cH^kg_+\cM$ are naturally equipped with an action of $\rW$ by automorphisms, we have $a_{\cH^kg_+\cM}=\cH^kg_+ a_\cM$ and $(\cH^kg_+\cM)^\ant=\cH^kg_+(\cM^\ant)$.

A similar result holds for $\QQ$-perverse sheaves on reduced complex analytic spaces and perverse cohomology sheaves of the direct image.
\end{proposition}

\begin{proof}
The first two points are clear by functoriality. We then have $\cH^kg_+\ker a_\cM\subset \ker a_{\cH^kg_+\cM}$ and $\cH^kg_+(\cM^\ant)\subset (\cH^kg_+\cM)^\ant$, and as the sum of both modules is equal to $\cH^kg_+\cM$, we get the third assertion.
\end{proof}

\subsection{Alternate Thom-Sebastiani for perverse sheaves}

Let $X$ be a reduced complex analytic space and let $f:X\to\CC$ be a holomorphic function. Let $\cF$ be a perverse sheaf of $\QQ$-vector spaces on $X$. Consider the $r$-fold product $X^r=X\times\cdots\times X$ with the function $f^{\oplus r}\defin f\oplus\cdots\oplus f:X^r\to\CC$ defined by
\[
f^{\oplus r}(x_1,\dots,x_r)=f(x_1)+\cdots+f(x_r),
\]
and the perverse sheaf $\cF^{\boxtimes r}\defin \cF\boxtimes\cdots\boxtimes\cF$.

Denote by $X^{(r)}$ the quotient space\footnote{$X^{(r)}$ is usually denoted by $\Sym^rX$, but we do not use the latter notation to avoid any confusion with \S\ref{subsec:alternate}.} of $X^r$ by the natural action of the symmetric group $\gS_r$ and let $\rho:X^r\to X^{(r)}$ be the projection. The space $X^{(r)}$ is a reduced analytic space (usually singular along the image of the diagonals, even if $X$ is smooth). The function $f^{\oplus r}$, being invariant under $\gS_r$, defines a holomorphic function $f^{(\oplus r)}:X^{(r)}\to\CC$ such that $f^{\oplus r}=f^{(\oplus r)}\circ\rho$.

The complex $\cG\defin\bR\rho_*\cF^{\boxtimes r}$ is a perverse sheaf (as $\rho$ is finite) and comes equipped with an action of $\gS_r$. We denote by $\cF^{\wedge r}=\cG^\ant$ its anti-invariant part (in the perverse category). If $D\cF$ denote the Verdier dual of $\cF$ on $X$, we have
\begin{equation}\label{eq:Dant}
\begin{split}
D(\cF^{\wedge r})=D(\cG^\ant)&\simeq (D\cG)^\ant\quad\text{according to Remark \ref{rem:antisym}},\\
& \simeq\big(\bR\rho_*(D\cF^{\boxtimes r})\big)^\ant\quad\text{as $\rho$ is finite},\\
&\simeq\big(\bR\rho_*((D\cF)^{\boxtimes r})\big)^\ant=(D\cF)^{\wedge r}.
\end{split}
\end{equation}

\subsubsection*{The case $\dim X=0$}
We assume that $X$ is a finite set of points. A $\QQ$-perverse sheaf $\cF$ on $X$ is then nothing but the data of a finite dimensional $\QQ$-vector space $F_x$ for each $x\in X$.
\begin{enumerate}
\item
If $X$ is reduced to a point $\{x\}$, and if we set $F=F_x$, then $X^{(r)}$ is reduced to a point and we have $\cF^{\wedge r}=\wedge^rF$.
\item
If $X$ is finite, we use the compatibility with the direct image $X\to\mathrm{pt}$ to see that $\Gamma(X^{(r)},\cF^{\wedge r})=\wedge^r(\oplus_{x\in X}F_x)$. If $x^{(r)}=\rho(x_1,\dots,x_r)$ is a point of $X^{(r)}$, the germ of $\cF^{\wedge r}$ at $x^{(r)}$ is the subspace of $\wedge^r(\oplus_{x\in X}F_x)$ generated by the $v_1\wedge\cdots\wedge v_r$, where $v_1\in F_{x_1},\dots,v_r\in F_{x_r}$, that we denote by $F_{x_1}\wedge\cdots\wedge F_{x_r}$.
\end{enumerate}

\begin{example}\label{ex:dim1}
Assume that $X$ is finite and $\dim F_x=1$ for any $x\in X$. Let $D\subset X^{(r)}$ be the image of the diagonals in $X^r$. Then $\cF^{\wedge r}_{|D}=0$. Indeed, if $x_1=x_2$ for instance, then $F_{x_1}=F_{x_2}$ and $F_{x_1}\wedge F_{x_2}=0$.
\end{example}

\subsubsection*{Restriction to a subset}
Let $i_Y:Y\hto X$ be the inclusion of a closed analytic subset. Let $\cF$ be a perverse sheaf on $X$. Assume that $i_Y^{-1}\cF$ is perverse up to a shift.

\begin{lemma}\label{lem:restr}
Under these assumptions, we have $(i_Y^{-1}\cF)^{\wedge r}=i_{Y^{(r)}}^{-1}(\cF^{\wedge r})$, where $i_{Y^{(r)}}$ is the natural inclusion $Y^{(r)}\hto X^{(r)}$.
\end{lemma}

\begin{proof}
Assume that $i_Y^{-1}\cF[k]$ is perverse, for some $k\in\ZZ$. Then $i_{Y^r}^{-1}\cF^{\boxtimes r}[kr]$ is perverse. On the other hand, we have $\cG_Y\defin\bR\varpi_*i_{Y^r}^{-1}\cF^{\boxtimes r}=i_{Y^{(r)}}^{-1}\bR\rho_*\cF^{\boxtimes r}=:i_{Y^{(r)}}^{-1}\cG$, as the diagram
\[
\xymatrix{
X^{r\vphantom{(r)}}\ar[r]^-{\rho}&X^{(r)}\\
Y^{r\vphantom{(r)}}\ar[r]_-{\varpi}\ar@{_{ (}->}[u]_{i_{Y^r}}&Y^{(r)}\ar@{_{ (}->}[u]_{i_{Y^{(r)}}}
}
\]
is cartesian. Then $i_{Y^{(r)}}^{-1}\cG[kr]$ is perverse. The decomposition $\cG=\ker a_\cG\oplus\cG^\ant$ induces a similar decomposition after applying the functor $i_{Y^{(r)}}^{-1}[kr]$, and we conclude as in Proposition \ref{prop:imdirant}.
\end{proof}

\subsubsection*{Fibre of $\cF^{\wedge r}$}
Assume that, up to a fixed shift, the restriction at $x_1,\dots,x_r$ of the perverse sheaf $\cF$ is a sheaf (in the following, we forget about the shift, which applies uniformly to all the sheaves that we consider). Applying Lemma \ref{lem:restr} to $Y=\{x_1,\dots,x_r\}$ and the case $\dim X=0$, we get
\[
\cF^{\wedge r}_{\rho(x_1,\dots,x_r)}\simeq\cF_{x_1}\wedge\cdots\wedge\cF_{x_r} \subset \wedge^r[\cF_{x_1}\oplus\cdots\oplus\cF_{x_r}].
\]
If for instance $\cF$ is a shifted local system of rank one, we can apply Example \ref{ex:dim1} to obtain that $\cF^{\wedge r}$ is a shifted local system of rank one on the complement of the image of the diagonals in $X^{(r)}$ and is zero on this image.

\subsubsection*{Application}
Let $X$ be a complex manifold, let $\pQQ_X=\QQ_X[\dim X]$ be the constant sheaf shifted by $\dim X$ (this is a perverse sheaf). Let us describe the perverse sheaf $\pQQ_X^{\wedge r}$. Denote by $D\subset X^{(r)}$ the image by $\rho$ of the diagonals of $X^r$ and by $V$ the open set $X^{(r)}\moins D$ (notice that it is smooth). Let $\delta:V\hto X^{(r)}$ denote the open inclusion.

As $\rho_*\QQ_{X^r}$ is a sheaf equipped with an action of $\gS_r$, we can also consider the anti-invariant subsheaf $(\rho_*\QQ_{X^r})^\ant$ (in the sense of sheaf theory). We denote it by~$\QQ_X^{\wedge r}$.

\begin{proposition}\label{prop:dim1lisse}\mbox{}
\begin{enumerate}
\item\label{prop:dim1lisse0}
We have $\pQQ_X^{\wedge r}=\QQ_X^{\wedge r}[r\dim X]$.
\item\label{prop:dim1lisse1}
The sheaf $\pQQ_{X|V}^{\wedge r}$ is a rank-one local system on $V$ shifted by $r\dim X$.

\item\label{prop:dim1lisse2}
With respect to Poincar\'e-Verdier duality, the perverse sheaf $\pQQ_X^{\wedge r}$ is self-dual.

\item\label{prop:dim1lisse3}
We have $\pQQ_X^{\wedge r}=\delta_!\delta^{-1}\pQQ_X^{\wedge r}=\bR\delta_*\delta^{-1}\pQQ_X^{\wedge r}$.
\end{enumerate}
\end{proposition}

\begin{proof}
Let us compute the germ of $\pQQ_X^{\wedge r}$ at some point $x^{(r)}$ of $X^{(r)}$. Denote by $Y=\vert x^{(r)}\vert\subset X$ the support of $x^{(r)}$. This is a finite set of points. We can apply Lemma \ref{lem:restr} to it, and then we can apply Example \ref{ex:dim1}. This shows \eqref{prop:dim1lisse1} and the first equality in \eqref{prop:dim1lisse3}. The second equality in \eqref{prop:dim1lisse3} is a consequence of the first one and of Poincar\'e duality \eqref{prop:dim1lisse2}.

Poincar\'e duality \eqref{prop:dim1lisse2} follows from \eqref{eq:Dant} and the self-duality of $\pQQ_X$.

Except from Poincar\'e duality, similar arguments can be applied to $\QQ_X^{\wedge r}$, showing that $\QQ_X^{\wedge r}=\delta_!\delta^{-1}\QQ_X^{\wedge r}$. It is moreover clear that $\delta^{-1}\pQQ_X^{\wedge r}=\delta^{-1}\QQ_X^{\wedge r}[r\dim X]$. This completes the proof of \eqref{prop:dim1lisse0}
\end{proof}

\begin{remark}
The complex $\pQQ_X^{\wedge r}$ is thus also equal to the intermediate extension $\delta_{!*}\delta^{-1}\pQQ_X^{\wedge r}$ (\ie the intersection complex attached to the rank-one shifted local system $\delta^{-1}\pQQ_X^{\wedge r}$).
\end{remark}

\begin{example}
Assume that $X=\Afu$. Then the space $(\Afu)^{(r)}$ is an affine space isomorphic to $\AA^{\!r}$. Let $D\subset \AA^{\!r}$ be the discriminant hypersurface and denote by $\delta:\AA^{\!r}\moins D\hto \AA^{\!r}$ the open inclusion.

The sheaf $\delta^{-1}\QQ_{\Afu}^{\wedge r}$ is a local system of rank one on $\AA^{\!r}\moins D$ with monodromy equal to $-\id$ locally around the smooth part of $D$ and we have
\[
\QQ_{\Afu}^{\wedge r}=\delta_!\delta^{-1}\QQ_{\Afu}^{\wedge r}=\bR\delta_*\delta^{-1}\QQ_{\Afu}^{\wedge r}.
\]
\end{example}

\begin{example}[Vanishing cycles]\label{exam:vanishcycles}
Let us come back to the case of a general perverse sheaf $\cF$ on~$X$. Denote by $\phip$ the functor of vanishing cycles shifted by $-1$ (see \eg \cite{Dimca04}). This is an exact functor on $\Perv(X)$. Let $C\subset X$ be the set of critical points of $f$ with respect to $\cF$: by definition $x^o\in C$ iff the germ at $x^o$ of the perverse sheaf $\phip_{f-f(x^o)}\cF$ is non-zero. Let us assume that $C$ is finite. Then, for any $x^o\in C$, the germ of $\phip_{f-f(x^o)}\cF$ at $x^o$ is the direct image by the embedding $\{x^o\}\hto X$ of a finite dimensional vector space $E_{x^o}$ (vanishing cycles of $f$ at $x^o$).

Let $x_1,\dots,x_r\in X$. According to the Thom-Sebastiani isomorphism for perverse sheaves with monodromy (see \cite{Massey01,Schurmann03}),
\[
\phip_{(f-f(x_1))\oplus\cdots\oplus(f-f(x_r))}(\cF^{\boxtimes r})_{(x_1,\dots,x_r)}\isom\phip_{(f-f(x_1))}\cF_{x_1}\boxtimes\cdots\boxtimes\phip_{(f-f(x_r))}\cF_{x_r}.
\]
It follows that $C^r$ is the set of critical points of $f^{\oplus r}$ (hence is finite) and that, if we denote by $\phip_{f,\tot}\cF$ the direct sum $\oplus_{x\in C}\phip_{f-f(x)}\cF$, then $\phip_{f^{\oplus r},\tot}(\cF^{\boxtimes r})\isom(\phip_{f,\tot}\cF)^{\boxtimes r}$.

For any complex number $c$, $\bR\rho_*\phip_{f^{\oplus r}-c}=\phip_{f^{(\oplus r)}-c}\bR\rho_*$ and, as $\phip_{f^{(\oplus r)}-c}$ is an exact functor on $\Perv(X^{(r)})$, it commutes with ${}^\ant$. It follows that the set of critical points of $f^{(\oplus r)}$ with respect to $\cF^{\wedge r}$ is contained in $C^{(r)}$, which is also finite. We then get
\begin{equation}\label{eq:phitot}
\phip_{f^{(\oplus r)},\tot}(\cF^{\wedge r})\isom(\phip_{f,\tot}\cF)^{\wedge r}
\end{equation}
More precisely, let us set $C=\{x_i\mid i\in\nobreak I\}$ and let us choose some total order on $I$. The set $C^{(r)}$ consists of the points $x_{i_1,\dots,i_r}=\rho(x_{i_1},\dots,x_{i_r})$ with $i_1\leq\cdots\leq i_r$. Let us set $E_i=E_{x_i}$. The critical value of $f^{(\oplus r)}$ at $x_{i_1,\dots,i_r}$ is $f^{(\oplus r)}(x_{i_1,\dots,i_r})=f(x_{i_1})+\cdots+f(x_{i_r})$. Then, according to the previous results, the space $(\phip_{f^{(\oplus r)}-f^{(\oplus r)}(x_{i_1,\dots,i_r})}(\cF^{\wedge r}))_{x_{i_1,\dots,i_r}}$ of vanishing cycles of $f^{(\oplus r)}$ at $x_{i_1,\dots,i_r}$ relatively to $\cF^{\wedge r}$ is the $(i_1,\dots,i_r)$-component of the alternate product $\wedge^r(\oplus_{i\in I}E_i)$.

Assume that all critical points of $f$ are simple (\ie $\dim E_i=1$ for any $i\in I$). Then the space of vanishing cycles of $f^{(\oplus r)}$ at $x_{i_1,\dots,i_r}$ relatively to $\cF^{\wedge r}$ vanishes as soon as two indices $i_a$ and $i_b$ (with $a\neq b$) coincide.
\end{example}

\subsection{Alternate Thom-Sebastiani for cohomologically tame functions}\label{subsec:nonchar}
Let $f:U\to\Afu$ be a regular function on a smooth affine algebraic variety~$U$.\footnote{When considering perverse sheaves, we implicitly use the underlying analytic objects.} Assume that there exists an algebraic variety $X$ in which $U$ is Zariski dense and a projective morphism $F:X\to\Afu$ inducing $f$ on $U$, such that, denoting by $j$ the inclusion $U\hto X$, for any $c\in\Afu$ the complex $\phip_{F-c}\bR j_*\pQQ_U$ is supported on a finite set of points in $U$ (that is, in $f^{-1}(c)$). We then say that $f$ is cohomologically tame with respect to the constant sheaf $\pQQ_U$.

In the remaining of this section we will assume that $f$ is cohomologically tame with respect to the constant sheaf.

\begin{lemma}\label{lem:noncharW}
If $f$ is cohomologically tame with respect to $\pQQ_U$, then $f^{(\oplus r)}$ is so with respect to the perverse sheaf $\pQQ_U^{\wedge r}$.
\end{lemma}

\begin{proof}
As $\pQQ_U^{\wedge r}$ is a direct summand of $\bR\rho_*\pQQ_{U^r}$, it is enough to prove the assertion for the latter perverse sheaf. Using the isomorphism $\bR\rho_*\phip_{F^{\oplus r}-c}=\phip_{F^{(\oplus r)}-c}\bR\rho_*$, we are reduced to proving the assertion for the Thom-Sebastiani sum $f^{\oplus r}$ with respect to $\pQQ_{U^r}$ and the partial compactification $U^r\Hto{j^r} X^r\To{F^{\oplus r}}\CC$. The result follows then from the Thom-Sebastiani theorem of \cite{Massey01,Schurmann03}, applied to $\bR j^r_*\pQQ_{U^r}$ and $F^{\oplus r}$.
\end{proof}

\begin{remark}
Let us assume that the critical points of $f$ are simple. It follows then from Example \ref{exam:vanishcycles} that the restriction $f^{(\oplus r)}$ to the open set $V\subset U^{(r)}$ is cohomologically tame with respect to the rank-one local system $\QQ_{U|V}^{\wedge r}$, and its critical points are simple. However, even if $f$ has distinct critical values, this may not remain true for $f^{(\oplus r)}$. Let us also notice that $V$ is smooth but usually not affine.
\end{remark}

\subsubsection*{The alternate Gauss-Manin system}
Let us recall how the Gauss-Manin system $G_f$ is defined from $f$. One first defines the differential system $M_f$ on the affine line with coordinate $t$ by setting
\[
M_f=\Omega^n(U)[\partial_t]\Big/(d-\partial_t df\wedge)\Omega^{n-1}(U)[\partial_t]\qquad (n=\dim U).
\]
This is a finite $\CC[t]\langle\partial_t\rangle$-module with regular singularities. From the point of view of $\cD$-modules, it is the direct image $\cH^0f_+\cO_U$, where we regard $\cO_U$ as a left module on the sheaf of differential operators $\cD_U$. As a consequence, the analytic de~Rham complex of $M_f$ is the $0$-th perverse cohomology of the direct image $\bR f_*\pCC_U$, where we denote as above by $\pCC_U$ the constant sheaf shifted by $\dim U$.

On the other hand, $G_f$ is interpreted as the localized algebraic Laplace transform of $M_f$: we set $\hb=\partial_t$, $\partial_\hb=-t$, so
\[
G_f=\Omega^n(U)[\hb,\hb^{-1}]\Big/(d-\hb df\wedge)\Omega^{n-1}(U)[\hb,\hb^{-1}].
\]
It is a free $\CC[\hb,\hb^{-1}]$-module, whose rank is equal to the dimension of the relative cohomology $H^n(U,f^{-1}(c))$ for a generic $c\in\Afu$. The Brieskorn lattice
\[
G_{f,0}\defin\Omega^n(U)[\hbm]\Big/(\hbm d-df\wedge)\Omega^{n-1}(U)[\hbm]
\]
is a free $\CC[\hbm]$-submodule of the same rank, where we set $\hbm=\hb^{-1}$ (\cf \cite{Bibi96bb}).

Let us begin with the tensor product:

\begin{lemma}[\cite{N-S97}]\label{lem:tensorTS}
The $r$-fold tensor product $\otimes^r_{\CC[\hbm]}G_{f,0}$ is isomorphic to the Bries\-korn lattice system of the $r$-fold Thom-Sebastiani sum $f^{\oplus r}:U^r\to\Afu$, where we set $f^{\oplus r}(u^{(1)},\dots,u^{(r)})=f(u^{(1)})+\cdots+f(u^{(r)})$.
\end{lemma}

\begin{proof}
Let us recall the proof. By an easy induction on $r$, it is enough to prove the result for the tensor product corresponding to cohomologically tame functions $f:U\to\Afu$ and $g:U'\to\Afu$. We consider the complex $(\Omega^{n+\cbbullet}(U)[\hbm],\hbm d-df\wedge)$. As $\hbm d-df\wedge$ is the twisted differential $e^{f/\hbm}\circ d\circ e^{-f/\hbm}$, we can write this complex as $(\Omega^{n+\cbbullet}(U)[\hbm]e^{-f/\hbm},d)$, where $e^{-f/\hbm}$ is now a symbol denoting the twist of the differential (to follow the definition of a shifted complex, we should use the differential $(-1)^nd$, but it is of no use here).

We have a natural morphism of complexes
\begin{multline*}
\tag{$*$}
(\Omega^{n+\cbbullet}(U)[\hbm]e^{-f/\hbm},d)\otimes_{\CC[\hbm]} (\Omega^{m-\cbbullet}(U')[\hbm]e^{-g/\hbm},d)\\
\to(\Omega^{n+m-\cbbullet}(U\times U')[\hbm]e^{-(f\oplus g)/\hbm},d).
\end{multline*}
It induces a surjective morphism of the corresponding $H^0$ as $\CC[\hbm]$-modules, since
\begin{multline*}
H^0\Big(\Omega^{n+m-\cbbullet}(U\times U')[\hbm]e^{-(f\oplus g)/\hbm},d\Big)\\
=\Omega^{n+m}(U\times U')[\hbm]\Big/(\hbm d-d(f\oplus g)\wedge)\Omega^{n+m-1}(U\times U')[\hbm]\
=G_{f\oplus g,0}.
\end{multline*}
As we have seen above (\cf also \cite[\S2]{N-S97}), $f\oplus g$ is cohomologically tame, hence $G_{f\oplus g,0}$ is a free $\CC[\hbm]$-module of finite rank (\cf \cite{Bibi96bb}).
On the other hand, using that $G_{f,0}$ and $G_{g,0}$ are free $\CC[\hbm]$-modules of finite rank (as a consequence of cohomological tameness), we identify the $H^0$ of the left-hand term in $(*)$ to $G_{f,0}\otimes_{\CC[\hbm]}G_{g,0}$, which is also free. We thus have a surjective morphism $G_{f,0}\otimes_{\CC[\hbm]}G_{g,0}\to G_{f\oplus g,0}$ of free $\CC[\hbm]$-modules. Moreover, a simple computation shows that their rank is the same. Therefore, this morphism is an isomorphism.
\end{proof}

\begin{example}\label{exam:torus}
Let us give the explicit description of the action of the symmetric group $\gS_r$ on $G_{f^{\oplus r},0}$ coming from the isomorphism of Lemma \ref{lem:tensorTS}, when $U$ is the torus $(\CC^*)^n$ with coordinates $u_1,\dots,u_n$ and volume form $\vol=\frac{du_1}{u_1}\wedge\cdots\wedge\frac{du_n}{u_n}$. For $\omega\in\Omega^n(U)[\hbm]$, we write $\omega=\varphi(u)\vol$ with $\varphi(u)\in\cO(U)[\hbm]$. Then $\omega_1\otimes\cdots\otimes\omega_r$ is sent to
\[
\varphi_1(u^{(1)})\cdots\varphi_r(u^{(r)})\vol_1\wedge\cdots\wedge\vol_r,
\]
and $\rw(\omega_1\otimes\cdots\otimes\omega_r)=\omega_{\rw(1)}\otimes\cdots\otimes\omega_{\rw(r)}$ is sent to
\[
\varphi_{\rw(1)}(u^{(1)})\cdots\varphi_{\rw(r)}(u^{(r)})\vol_1\wedge\cdots\wedge\vol_r.
\]
Therefore, after dividing by $\vol_1\wedge\cdots\wedge\vol_r$, the action of $\gS_r$ on $G_{f^{\oplus r},0}$ amounts to the usual action induced by that on $\cO(U^r)$.
\end{example}

We are now interested in the alternate product $\wedge^r_{\CC[\hbm]}G_{f,0}$, that is, the antisymmetric submodule of $G_{f^{\oplus r},0}$. In the situation of Example \ref{exam:torus}, it is isomorphic to $\cO(U^r)^\ant[\hbm]$. We will now give another interpretation of this antisymmetric submodule, or at least of the corresponding submodule of the Gauss-Manin system $G_{f^{\oplus r}}$, as the Gauss-Manin system attached to the morphism $f^{(\oplus r)}$ induced by $f^{\oplus r}$ on the quotient variety $U^{(r)}=\Ur/\gS_r$, with respect to the perverse sheaf $\pCC_U^{\wedge r}=\CC_U^{\wedge r}[r\dim U]$ described by Proposition \ref{prop:dim1lisse}. Let us note that the quotient variety $U^{(r)}$ is affine, but usually singular.

Let us consider the $0$-th perverse cohomology $\pcH^0\bR f^{(\oplus r)}_*\pCC_U^{\wedge r}$. This is a perverse sheaf on $\Afu$, which corresponds to a unique (up to isomorphism) regular $\CC[t]\langle\partial_t\rangle$-module that we denote by $M_f^{\wedge r}$.

\begin{proposition}\label{prop:laplacewedge}
The localized algebraic Laplace transform of $M_f^{\wedge r}$ is isomorphic to $\wedge^rG_f$.
\end{proposition}

\begin{proof}[Sketch of proof]
Let us choose an embedding $U^{(r)}\hto \cU$ into a smooth affine variety and let us still denote by $\rho$ the finite morphism $U^r\to \cU\times\Afu$ obtained by composing $\rho:U^r\to U^{(r)}$ with the graph embedding $\iota:U^{(r)}\hto \cU\times\Afu$ of $f^{(\oplus r)}$.

We will work with right $\cD$-modules and we denote by $\omega_U$ the right $\cD_U$-module $\Omega^{\dim U}_U$. The $\cD_{U^r}$-module ${\boxtimes}^r\omega_U$ is $\gS_r$-equivariant, so $\rho_+({\boxtimes}^r\omega_U)$ has an action of~$\gS_r$. Taking Spencer complexes (which plays the role of the de~Rham complex for right $\cD$-modules), we have an isomorphism $\Sp_{\cU\times\Afu}^\cbbullet(\rho_+({\boxtimes}^r\omega_U))=\rho_*\pCC_{U^r}$ which is compatible with the $\gS_r$-action. Therefore,
\[
\Sp_{\cU\times\Afu}^\cbbullet(\omega_U^{\wedge r})=\iota_*\pCC_U^{\wedge r}.
\]
Using the compatibility with direct images we find
\[
\bR f^{(\oplus r)}_*\pCC_U^{\wedge r}=\Sp_{\Afu}^\cbbullet\big(f^{(\oplus r)}_+(\omega_U^{\wedge r})\big)=\Sp_{\Afu}^\cbbullet\big((f^{\oplus r}_+({\boxtimes}^r\omega_U))^\ant)\big).
\]
Therefore,
\[
\Sp_{\Afu}^\cbbullet M_f^{\wedge r}\defin\pcH^0\bR f^{(\oplus r)}_*\pCC_U^{\wedge r}\simeq\Sp_{\Afu}^\cbbullet\big(\cH^0(f^{\oplus r}_+({\boxtimes}^r\omega_U))^\ant\big).
\]
If we set $M_f^{\oplus r}=\cH^0(f^{\oplus r}_+({\boxtimes}^r\omega_U))$, we thus have $M_f^{\wedge r}\simeq(M_f^{\oplus r})^\ant$. On the other hand, by definition, the localized Laplace transform of $M_f^{\oplus r}$, with its $\gS_r$-action, is isomorphic, by Lemma \ref{lem:tensorTS}, to $\otimes^r_{\CC[\hb,\hb^{-1}]}G_f$ with the natural action of~$\gS_r$.
\end{proof}

\section{Alternate Thom-Sebastiani and Frobenius manifolds}\label{sec:geom}

In this section, we consider a function $f:U\to\Afu$ satisfying the assumptions of \S\ref{subsec:nonchar}.

\subsection{The canonical pre-Saito structure}\label{subsec:canSaito}
We denote by $\ov G_f$ the $\CC$-vector space $G_f$ on which the action of $\CC[\hb,\hb^{-1}]\langle\partial_\hb\rangle$ is modified by a sign: we set, for any $g\in G_f$, $\hb\cdot g=-\hb g$ and $\partial_\hb\cdot g=-\partial_\hb g$. In other words,
\[
\ov G_f=\Omega^n(U)[\hb,\hb^{-1}]\Big/(d+\hb df\wedge)\Omega^{n-1}(U)[\hb,\hb^{-1}]
\]
equipped with the action of $\partial_\hb$ defined by $\partial_\hb[\omega]=[f\omega]$ for any $\omega\in\Omega^n(U)$. The Brieskorn lattice $\ov G_{f,0}$ is defined similarly.

We will use the following two results (\cf \cite{Bibi96bb}):
\begin{enumerate}
\item
Poincar\'e duality for the morphism $f$ induces a \emph{canonical} nondegenerate $(-1)^n$-Hermitian sesquilinear pairing
\[
S_f:G_f\otimes_{\CC[\hb,\hb^{-1}]}\ov G_f\to\CC[\hb,\hb^{-1}]
\]
which is compatible with the action of $\partial_\hb$ (that is, $\partial_\hb(S_f(g',g''))=S_f(\partial_\hb g',g'')-S_f(g',\partial_\hb g'')$). This pairing induces a perfect pairing
\[
S_f:G_{f,0}\otimes_{\CC[\hbm]}\ov G_{f,0}\to\hbm^n\CC[\hbm].
\]

\item
The limit mixed Hodge structure on $\lim_{c\to\infty}H^n(U,f^{-1}(c),\QQ)$ enables one to produce, through a construction due to M.~Saito \cite{MSaito89}, a \emph{canonical} $\CC$-vector space $E^o_f$ of $G_{f,0}$, such that $G_{f,0}=\CC[\hbm]\otimes_\CC E^o_f$, and in which $\partial_\hb=-\hbm^2\partial_\hbm$ takes the form
\[
-R_0^o+\frac{R_\infty}{\hb},
\]
where $R_0^o$ and $R_\infty$ are two endomorphisms of $E^o_f$. Moreover, restricting $S_f$ to $E^o_f\otimes_\CC E^o_f$ gives a symmetric nondegenerate pairing
\[
g^o:E^o_f\otimes_\CC E^o_f\to\CC.
\]
Lastly, on the one hand, $R_\infty$ is semisimple and its spectrum is the opposite of the \emph{spectrum at infinity} of $f$; if $R_\infty^*$ denotes its $g^o$-adjoint, we have $R_\infty+R_\infty^*=-n\id$. On the other hand, through the isomorphism $E^o_f\to G_{f,0}/\hbm G_{f,0}=\Omega^n(U)/df\wedge\Omega^{n-1}(U)$, $R_0^o$ corresponds to the endomorphism induced by the multiplication by $f$ on $\Omega^n(U)$ and satisfies $R_0^{o*}=R_0^o$. Its eigenvalues, counted with multiplicity, are the critical values of $f$ counted with multiplicity.
\end{enumerate}

In other words, to any such function $f$ is associated (mainly using Hodge theory related to it) a \emph{canonical} pre-Saito structure $(E^o_f,R_\infty,R_0^o,g^o)$ of weight \hbox{$w=\dim U=n$} with base manifold reduced to a point.

\subsection{The trivial deformation}
We now show that the trivial deformation, as constructed in Example \ref{exam:deformation}, of the pre-Saito structure $(E^o_f,R_\infty,R_0^o,g^o)$ defined above can be obtained from a deformation of $f$ itself.

Let $\CC^*$ be the one-dimensional torus with coordinate $\lambda$. Later, we will consider the analytic uniformization $\lambda=e^x$ to be compatible with Example \ref{exam:deformation}. For $f$ as in \S\ref{subsec:nonchar}, we consider the unfolding
\[
F:U\times\CC^*\to\Afu,\qquad (u,\lambda)\mto \lambda f(u).
\]
The Gauss-Manin system $G_F$ of $F$ is a one-parameter deformation of that of $f$. We set (still denoting by $d$ the differential with respect to the $U$-variables only)
\[
G_{F,0}=\Omega^n(U)[\lambda,\lambda^{-1},\hbm]\Big/(\hbm d-\lambda df\wedge)\Omega^{n-1} (U)[\lambda,\lambda^{-1},\hbm].
\]
The action of $\hbm^2\partial_\hbm$ is induced by the multiplication by $\lambda f$ on $\Omega^n(U)[\lambda,\lambda^{-1}]$ (and extended with the Leibniz rule). The action of $\zeta\partial_\lambda$ is induced by the multiplication by $-f$ on $\Omega^n(U)[\hbm]$ (and extended with the Leibniz rule).

Let us denote by $\pi$ the map $(\lambda,\hb)\mto \lambda\hb$ and by $\pi^*:\CC[\hb]\to\CC[\lambda,\lambda^{-1},\hb]$ or $\CC[\hbm]\to\CC[\lambda,\lambda^{-1},\hbm]$ the corresponding morphism of algebras, defined by $\hb\mto \lambda\hb$ and $\hbm\mto \lambda^{-1}\hbm$. Then $G_{F,0}=\pi^+G_{f,0}$, where $\pi^+$ means $\pi^*$ of the $\CC[\hbm]$-module and the natural lifting of the connection. Regarding $\CC[\lambda,\lambda^{-1},\hbm]$ as a $\CC[\hbm]$-module through $\pi^*$, we have $G_{F,0}=\CC[\lambda,\lambda^{-1},\hbm]\otimes_{\CC[\hbm]}G_{f,0}$ and
\[
\hbm^2\partial_\hbm(1\otimes g)=\lambda\otimes(\hbm^2\partial_\hbm g),\quad \hbm\partial_\lambda(1\otimes g)=-1\otimes(\hbm^2\partial_\hbm g).
\]
Using the space $E^o_f\subset G_{f,0}$ given by Hodge theory and M.~Saito's procedure for~$f$, we get the trivialization $G_{F,0}=\CC[\lambda,\lambda^{-1},\hbm]\otimes_\CC E^o_f$, and we get a pre-Saito structure by changing the trivialization as in Remark \ref{rem:deformation} (using here the variable~$\hbm$ instead of~$\hb$). From Remark \ref{rem:deformation} we obtain:

\begin{proposition}\label{prop:preFrobdef}
Let $(E^o_f,R_\infty,R_0^o,g^o)$ be the canonical pre-Saito structure of weight $n$ attached to $f$. Then, for any $x\in\Afu$, the canonical pre-Saito structure of weight $n$ attached to $e^xf$ is the fibre at $x$ of the trivial deformation of $(E^o_f,R_\infty,R_0^o,g^o)$ constructed in Example \ref{exam:deformation} (plus Example \ref{exam:deformationmetric} for the metric).\qed
\end{proposition}

\subsection{Frobenius manifold structure}
In order to obtain a Frobenius manifold, we need a pre-primitive homogeneous section $\omega^o$, canonically associated to the geometry. Such a section exists when $U$ is a torus, so we will only consider this case.

\begin{assumption}\label{assumpt:torus}
$U\simeq(\CC^*)^n$ is a torus with coordinates $u_1,\dots,u_n$ and $f:U\to\Afu$ is a Laurent polynomial such that:
\begin{enumerate}
\item
$f$ is convenient and nondegenerate with respect to its Newton polyhedron (\cf \cite{Kouchnirenko76}),
\item
the critical points of $f$ are simple and the critical values are distinct.
\end{enumerate}
\end{assumption}

As a consequence, the Jacobian algebra $\cO(U)/(\partial f)$ is finite dimensional, and the multiplication by $f$ induces on it a regular semisimple endomorphism, whose eigenvalues are the critical values of $f$. Moreover, $f$ is cohomologically tame with respect to the constant sheaf. We can apply to it the results indicated above (\cf \cite[\S4]{D-S02a}).

The class $\omega^o$ of the volume form $\frac{du_1}{u_1}\wedge\cdots\wedge\frac{du_n}{u_n}$ belongs to the canonically defined vector space $E^o_f$ and is homogeneous of degree $0$ with respect to $R_\infty$. Moreover, it is a cyclic vector for $R_0^o$. It is thus pre-primitive and homogeneous. Therefore, the data $(E^o_f,R_\infty,R_0^o,g^o,\omega^o)$ define a canonical Frobenius manifold structure on $\wt E^o_f$ of weight $n$. Let us note that any other coordinate system on the torus, obtained from $(u_1,\dots,u_n)$ by a monomial change of coordinates, leads to a new volume form equal to $\pm\omega^o$. According to Remark \ref{rem:sqrt}, the Frobenius structure does not depend on the choice of the coordinate system on the torus.

Let us now consider the $r$-fold alternate product. From Proposition \ref{prop:laplacewedge}, we get:

\begin{corollary}\label{cor:GMFrob}
The restriction to $\Afu$ of the differential system $(\bF^o,\nablab^o)$ on $\PP^1$ associated to the $r$-fold alternate product of the canonical Frobenius manifold attached to $f$ is the Gauss-Manin system $\wedge^rG_f$ of the pair $(f^{(\oplus r)},\pCC_U^{\wedge r})$ on $U^{(r)}$.\qed
\end{corollary}

\begin{example}
Let $f(u)=u_0+u_1+\cdots+u_n$, where we have set $u_0=1/(u_1\cdots u_n)$. The canonical pre-Saito structure $(E^o_f,R_\infty,R_0^o,g^o,\omega^o)$ is obtained in the following way (see for instance \cite{D-S02b} with all the weights set to one). The space $E^o_f$ is the $\CC$-vector space generated by $\omega^o_0=\omega^o,\omega^o_1,\dots,\omega^o_n$, where, for $k\geq1$, $\omega^o_k$ is the class of $u_0\cdots u_{k-1}\,\frac{du_1}{u_1}\wedge\cdots\wedge\frac{du_n}{u_n}$. In this basis, the matrices of $R_\infty,R_0^o,g^o$ are those of Example \ref{exam:proj}.\footnote{A similar computation can be done with weights, for instance by setting $u_0=1/(u_1^{w_1}\cdots u_n^{w_n})$ with $w_1,\dots,w_n\in\NN^*$ (\cf \cite{D-S02b,Mann07,C-C-L-T06}).}

From Theorem \ref{th:cfks} and Corollary \ref{cor:GMFrob} we conclude that the Gauss-Manin system of the pair $(e^{(r-1)i\pi/(n+1)}f^{(\oplus r)},\pCC_U^{\wedge r})$ can also be obtained from the Frobenius manifold attached to the quantum cohomology of the Grassmannian at its origin.
\end{example}

\begin{remark}
It would be desirable to give an interpretation of $\wedge E^o_f$ and of the metric induced by $\otimes^rg^o$ purely in terms of $(f^{(\oplus r)},\pCC_U^{\wedge r})$ (by using Hodge theory at $f^{(\oplus r)}=\infty$), so that the canonical process of \S\ref{subsec:canSaito} could be directly applied to $(f^{(\oplus r)},\pCC_U^{\wedge r})$.

On the other hand, it would also be desirable to define a suitable small deformation of $f^{(\oplus r)}$ which would be enough to recover the $r$-fold alternate product of the pre-Saito structure attached to~$f$. A natural choice would be the deformation induced by the deformation of $f^{\oplus r}$ by the elementary symmetric functions of the $f(u^{(i)})$ ($i=1,\dots,r$), but this deformation usually introduces new critical points, which would have to be eliminated in some way.
\end{remark}

\providecommand{\bysame}{\leavevmode\hbox to3em{\hrulefill}\thinspace}

\end{document}